\documentclass[11pt]{amsart}
\usepackage{amssymb,amsmath,epsfig,graphics,euscript}
\theoremstyle{plain}
\newtheorem{thrm}{Theorem}[section]
\newtheorem{lemma}[thrm]{Lemma}
\newtheorem{prop}[thrm]{Proposition}
\newtheorem{cor}[thrm]{Corollary}
\newtheorem{rmrk}[thrm]{Remark}
\newtheorem{dfn}[thrm]{Definition}

\numberwithin{equation}{section} \numberwithin{figure}{section}

\begin{document}
\newcommand{\SL}{\mathcal L^{1,p}(\Om)}
\newcommand{\Lp}{L^p(\Omega)}
\newcommand{\CO}{C^\infty_0(\Omega)}
\newcommand{\Rn}{\mathbb R^n}
\newcommand{\Rm}{\mathbb R^m}
\newcommand{\R}{\mathbb R}
\newcommand{\Om}{\Omega}
\newcommand{\Hn}{\mathbb H^n}
\newcommand{\HH}{\mathbb H^1}
\newcommand{\eps}{\epsilon}
\newcommand{\BVX}{BV_H(\Omega)}
\newcommand{\IO}{\int_\Omega}
\newcommand{\bG}{\boldsymbol{G}}
\newcommand{\bg}{\mathfrak g}
\newcommand{\p}{\partial}
\newcommand{\Xnu}{\overset{\rightarrow}{ X_\nu}}
\newcommand{\nuX}{\boldsymbol{\nu}^H}
\newcommand{\Up}{\boldsymbol{\mathcal Y}_H}
\newcommand{\n}{\boldsymbol \nu}
\newcommand{\sigmau}{\boldsymbol{\sigma}^u_H}
\newcommand{\di}{\nabla^{H,\mS}_i}
\newcommand{\duno}{\nabla^{H,\mS}_1}
\newcommand{\ddue}{\nabla^{H,\mS}_2}
\newcommand{\del}{\delta_H}
\newcommand{\nui}{\nu_{H,i}}
\newcommand{\nuj}{\nu_{H,j}}
\newcommand{\dej}{\delta_{H,j}}
\newcommand{\cx}{\boldsymbol{c}^\mathcal S}
\newcommand{\sx}{\sigma_H}
\newcommand{\lx}{\mathcal L_H}
\newcommand{\pb}{\overline p}
\newcommand{\qb}{\overline q}
\newcommand{\ob}{\overline \omega}
\newcommand{\nuu}{\boldsymbol \nu_{H,u}}
\newcommand{\nuv}{\boldsymbol \nu_{H,v}}
\newcommand{\Bl}{\Bigl|_{\lambda = 0}}
\newcommand{\mS}{\mathcal S}
\newcommand{\delh}{\Delta^H}
\newcommand{\delinf}{\Delta_{H,\infty}}
\newcommand{\nabh}{\nabla^H}
\newcommand{\delp}{\Delta_{H,p}}
\newcommand{\mO}{\mathcal O}
\newcommand{\delhs}{\Delta_{H,S}}
\newcommand{\lhs}{\hat{\Delta}_{H,S}}
\newcommand{\bN}{\boldsymbol{N}}
\newcommand{\bnu}{\boldsymbol \nu}
\newcommand{\la}{\lambda}
\newcommand{\nup}{({\boldsymbol{\nu}^H})^\perp}
\newcommand{\oX}{\overline X}
\newcommand{\oY}{\overline Y}
\newcommand{\ou}{\overline u}
\newcommand{\dels}{\nabla^{H,\mS}}
\newcommand{\delXY}{\nabla^{H,\mS}_X Y}
\newcommand{\delYX}{\nabla^{H,\mS}_Y X}
\newcommand{\hp}{\tilde{h}_0'(s)}
\newcommand{\uc}{u_{,ij}}
\newcommand{\sij}{\sum_{i,j=1}^m}
\newcommand{\sul}{\Delta^H}
\newcommand{\nor}{\boldsymbol \nu}
\newcommand{\fv}{\mathcal V^{H}_I(\mS;\mathcal X)}
\newcommand{\sv}{\mathcal V^{H}_{II}(\mS;\mathcal X)}

\newcommand{\Id}{I_\delta}
\newcommand{\Ri}{\mathbb{R} \times I}
\newcommand{\Rd}{\mathbb{R} \times \Id}

\newcommand{\til}{\tilde{\boldsymbol{\nu}}^H}
\newcommand{\tilp}{(\tilde{\boldsymbol{\nu}}^H)^\perp}
\newcommand{\tilL}{\tilde{\mathcal{L}}_z(r)}
\newcommand{\tilLg}{\tilde{\mathcal{L}}_{\gamma(s)}(r)}
\newcommand{\tilLp}{\tilde{\mathcal{L}}_z'(r)}
\newcommand{\tilLzero}{\tilde{\mathcal{L}}_z(0)}
\newcommand{\ra}{\rightarrow}
\renewcommand{\to}{\rightarrow}

\title[Instability of graphical strips and a positive answer, etc.] {Instability of graphical strips and a positive answer to the Bernstein problem in the
Heisenberg group $\HH$}

\author{D. Danielli}
\address{Department of Mathematics\\Purdue University \\
West Lafayette, IN 47907} \email[Donatella
Danielli]{danielli@math.purdue.edu}
\thanks{First author supported in part by NSF CAREER Grant, DMS-0239771}

\author{N. Garofalo}
\address{Department of Mathematics\\Purdue University \\
West Lafayette, IN 47907} \email[Nicola
Garofalo]{garofalo@math.purdue.edu}
\thanks{Second author supported in part by NSF Grant DMS-0300477}

\author{D. M. Nhieu}
\address{Department of Mathematics\\Georgetown University \\
Washington, DC 20057-1233} \email[Duy-Minh
Nhieu]{nhieu@math.georgetown.edu}

\author{S. D. Pauls}
\address{Department of Mathematics, Dartmouth College, Hanover, NH 03755}
\email{scott.pauls@dartmouth.edu}
\thanks{Fourth author supported in part by NSF Grant DMS-0306752}

%
%
 \subjclass{}
\date{\today}

\maketitle


\baselineskip 13.5pt


\section{\textbf{Introduction}}

\vskip 0.2in

One of the most celebrated problems in geometry and calculus of
variations is the Bernstein problem, which asserts that a $C^2$
minimal graph in $\R^3$ must necessarily be an affine plane.
Following an old tradition, here minimal means of vanishing mean
curvature. Bernstein \cite{Be} established this property in 1915.
Almost fifty years later a new insight of Fleming \cite{Fle}
sparked a major development in the geometric measure theory
which, through the celebrated works \cite{DG3}, \cite{Al},
\cite{Sim}, \cite{BDG} culminated in the following solution of the
Bernstein problem.

\begin{thrm}\label{T:classicalB}
Let $\mS = \{(x,u(x)) \in \R^{n+1}| x\in \Rn ,  x_{n+1} = u(x)\}$ be
a $C^2$ minimal graph in $\R^{n+1}$, i.e., let $u\in C^2(\Rn)$ be a
solution of the minimal surface equation
\begin{equation}\label{ms}
div\left(\frac{Du}{\sqrt{1 + |Du|^2}}\right)\ =\ 0\ ,
\end{equation}
in the whole space. If $n\leq 7$, then there exist $a\in \Rn$,
$\beta \in \R$ such that $u(x) = <a,x> + \beta$, i.e., $\mS$ must
be an affine hyperplane. If instead $n\geq 8$, then there exist
non affine (real analytic) functions on $\Rn$ which solve
\eqref{ms}.
\end{thrm}

The purpose of this paper is to study, in the simplest model of a
sub-Riemannian space, the three-dimensional Heisenberg group $\HH$,
the structure of $C^2$ minimal graphs with empty characteristic
locus and which, on every compact set, minimize the horizontal
perimeter. As a corollary of our results we obtain a positive answer
to a sub-Riemannian analogue of the Bernstein problem. From the
perspective of geometry the relevance of $\HH$ lies in the fact that
this Lie group constitutes the simplest
 prototype of a class of graded nilpotent Lie groups which arise
 as ``tangent spaces" in the Gromov-Hausdorff limit of Riemannian
 spaces, see \cite{Be}, \cite{Gro1},\cite{Gro2}, \cite{Mon}, \cite{CDPT}. Furthermore, $\Hn$ is an
 interesting model of a metric space with a non-trivial geometry.

The Bernstein problem in $\HH$ has recently received increasing
attention, see \cite{GP}, \cite{CHMY}, \cite{CH}, \cite{RR1},
\cite{RR2}, \cite{HP}, \cite{ASV}, \cite{DGN3}, \cite{GS},
\cite{CHY}, \cite{BSV}. While we refer to the discussion below, to
section \ref{S:prelim}, and to the cited references for a detailed
description of the relevant geometric setting, it seems appropriate
here to provide the reader with some historical perspective and a
brief overview of the main contributions of the present paper.

When approaching the sub-Riemannian Bernstein problem one is
confronted with the fact that there exists smooth entire minimal
graphs which are not affine. For instance, in $\HH$ the surface
$\mS$ defined by $t = xy/2$ is minimal, in the sense that its
$H$-mean curvature, defined in \eqref{hmc} below, vanishes
identically. This is of course in sharp contrast with Theorem
\ref{T:classicalB}. What lurks in the dark here are two aspects: (i)
On the one hand, the so-called characteristic locus of the surface
$\mS$, i.e., the collection of points at which the fiber of the
subbundle generated by the horizontal distribution coincides with
the tangent space of $\mS$: on this set the horizontal Gauss map
becomes singular (for instance, the surface $t = xy/2$ has non-empty
characteristic locus); (ii) The drastically different nature of the
relevant minimal surface equation, with respect to the classical
case. The classical minimal surface equation \eqref{ms} is a
quasi-linear elliptic equation, whereas in $\HH$, away from the
characteristic locus, one has a degenerate quasi-linear hyperbolic
equation (with one vanishing eigenvalue). Recall that a quasi-linear
second order equation $a u_{xx} + 2 b u_{xy} +c u_{yy}= d$, where
$a, b, c,d$ depend on $x,y,u,u_x,u_y$, is called elliptic (on the
solution $u$) if $ac - b^2>0$, hyperbolic is $ac-b^2<0$. For
instance, for a $C^2$ graph $\mS \subset \R^3$ of the type $t =
u(x,y)$, \eqref{ms} is equivalent to
\[ (1 + u_y^2) u_{xx} - 2 u_x u_y u_{xy} + (1 + u_x^2) u_{yy}\ =\ 0\
,
\]
for which one has $ac-b^2 = 1 + u_x^2 +u_y^2 >0$. On the other hand
the condition that $\mS$ be $H$-minimal becomes, away from the
characteristic locus $\Sigma(\mS) = \{(x,y,t)\in \mS \mid u_y -
\frac{x}{2} = 0\ ,\ u_x + \frac{y}{2} = 0\}$,
\[
\left(u_y - \frac{x}{2}\right)^2 u_{xx} - 2 \left(u_x + \frac{y}{2}
u_y\right)\left(u_y - \frac{x}{2}\right) u_{xy} + \left(u_x +
\frac{y}{2}\right)^2 u_{yy}\ =\ 0\ .
\]

Since in this case $ac-b^2\equiv 0$, we conclude that the equation
is degenerate hyperbolic. It is interesting to observe that since,
away from the characteristic locus, we can always locally
parameterize $\mS$ in one of the two forms \eqref{X1graph},
\eqref{X2graph} below, then by \eqref{H3} below, in terms of the
function $\phi(u,v)$ the $H$-minimality of $\mS$ is expressed by the
pde
\[
\phi_{uu} + 2 \phi \phi_{uv} +\phi^2 \phi_{vv}\ =\ -\ \phi_v (\phi
\phi_v + 2 \phi_u)\ ,
\]
which again is degenerate hyperbolic.

The above example $t = xy/2$ is not isolated since a basic result
first proved in \cite{CHMY}, and also independently  with a
different proof in \cite{GP}, shows in particular that every
$H$-minimal entire graph over the horizontal plane $t=0$ must have
non-empty characteristic locus, see Theorem \ref{T:xygraphs} below.
These considerations suggest that, in order to have a reasonably
behaved horizontal Gauss map, in the sub-Riemannian Bernstein
problem one should impose the restriction that the horizontal bundle
$H\HH$ be transversal to the tangent bundle $T\mS$ at every point of
the surface, i.e., $\mS$ should have empty characteristic locus.
This, in turn, immediately imposes a restriction on the type of
planes which are appropriate in the sub-Riemannian Bernstein
problem. Since every plane $ax + by + ct = \gamma$, for which $c\neq
0$, possesses the isolated characteristic point
$(-2b/c,2a/c,\gamma/c))$, it is clear that we want to confine the
attention to the non-characteristic vertical planes
\begin{equation}\label{vp}
\tilde{P}_0\ =\ \{(x,y,t)\in \HH \mid ax
+ by = \gamma\}\ .
\end{equation}

The appropriateness of these planes is also confirmed by the
fundamental Rademacher-Stepanov type theorem of Pansu \cite{Pa}.
Specialized to the present setting, the latter states that if $F:\HH
\to \R$ is a Lipschitz map with respect to the Carnot-Carath\'eodory
distance associated with the subbundle $H\HH$, then $F$ is Pansu
differentiable at a.e. point $g=(x,y,t)\in \HH$ (w.r.t. Lebesgue
measure), and the Pansu differential is given by $F_*(x,y,t) = a x +
b y$, for some $a,b\in \R$. This result underscores the special role
of the vertical planes \eqref{vp} in sub-Riemannian geometric
measure theory. A closely related remarkable fact, discovered in
\cite{FSS1}, is that the blow-up \`a la De Giorgi of a set with
locally finite $H$-perimeter at a point of its reduced boundary is
again a plane such as \eqref{vp}. It is time to introduce a basic
definition.

\begin{dfn}\label{D:graph}
We say that a surface $\mS\subset \HH$ of class $C^2$ is an
\emph{entire graph} if there exists a plane $P$, having equation $ax
+ b y + ct =\gamma$, such that for every point $g_0\in P$ the
straight line passing through $g_0$ and parallel to the Euclidean
normal $\bN_e = (a,b,c)$ to $P$, $L(g_0) = \{g_0 + s \bN_e\mid s\in
\R\}$, intersects the surface $\mS$ in exactly one point.
\end{dfn}

The above considerations suggest the natural conjecture that if
$\mS\subset \HH$ is an entire $H$-minimal graph, with empty
characteristic locus, then $\mS$ should be a vertical plane $\tilde
P_0$ such as \eqref{vp}. Since, as we have mentioned, $H$-minimal is
intended in the sense that the horizontal mean curvature $\mathcal
H$ vanishes identically as a continuous function on $\mS$, it is
worth stressing that, thanks to the first variation formula in
Theorem \ref{T:variations} below, a $C^2$ surface with empty
characteristic locus is $H$-minimal if and only if it is a critical
point of the $H$-perimeter given by \eqref{per} below. However, for
$\HH$ the situation is very different than in Euclidean space. In
fact, in \cite{GP} the second and the fourth named authors
discovered the following counterexample to such a plausible
sub-Riemannian version of the Bernstein problem. The real analytic
surface
\begin{equation}\label{ce}
S\ =\ \{(x,y,t)\in \mathbb H^1 \mid x =  y \tan \tanh(t)\}
\end{equation}
is an entire $H$-minimal graph, with empty characteristic locus,
over the coordinate $(y,t)$-plane in $\mathbb H^1$. This example
seems to cast a dim light over the sub-Riemannian Bernstein problem.

There is however a deeper aspect of the problem which in the
classical case is confined to the background, but which, due to the
diverse nature in the sub-Riemannian setting of the relevant area
functional, the horizontal perimeter, might be playing an important
role. What could be happening, in fact, is that $H$-minimal surfaces
such as \eqref{ce} are unstable, in the sense that they are only
critical points, but not local minimizers of the $H$-perimeter. This
phenomenon, which goes back to the classical findings of Bernoulli,
see e.g. \cite{Ca}, has of course no counterpart in the Bernstein
problem in flat space, since, thanks to the convexity of the area
functional $A(\mS) = \int_\Om \sqrt{1 + |Du|^2} dx$, stability is
automatic for a solution to \eqref{ms}, see e.g. \cite{CM}. On the
other hand, stability had played an important role in the work on
complete minimal surfaces in $3$-manifolds by Fischer-Colbrie and
Schoen \cite{FCS}, who generalized Bernstein's theorem and proved,
in particular, that the only complete stable oriented minimal
surfaces in $\R^3$ are the planes. This latter result was also
independently obtained in \cite{DCP}.

For the functional which expresses the horizontal perimeter the
convexity fails. To see this negative phenomenon consider for
instance the situation in which $\mS$ is parameterized by
\eqref{X1graph}, then it was proved in \cite{ASV} that the
$H$-perimeter of $\mS$ is expressed by the functional
\[ P_H(\mS)\ =\ \mathcal P(\phi)\ =\ \int_\Om \sqrt{1 + \mathcal B_\phi(\phi)^2}\ du dv\
,
\] where $\mathcal B_\phi(\phi)$ denotes the nonlinear inviscid
Burger operator acting on $\phi$, see \eqref{burger}. If we set $\xi
= \phi_u$, $\eta=\phi_v$, then the integrand of the above functional
is given by
\[
F(\phi,\xi,\eta)\ =\ \sqrt{1 + (\xi + \phi \eta)^2}\ .
\]

The Hessian of $F$ is given by
\[
\nabla^2 F\ =\ F^{-3} \ \begin{pmatrix}\eta^2 &  \eta  & \phi \eta +
(\xi + \phi \eta) F^2
\\
\eta & 1 & \phi
\\
\phi \eta + (\xi + \phi \eta) F^2 & \phi & \phi^2
\end{pmatrix}\ ,
\]
and therefore $det\ \nabla^2 F = - (\xi + \phi \eta)^2 F^{-5}$,
which shows that $F$ is not convex.

These considerations suggest that the above conjecture could be
repaired as follows: \emph{In $\HH$ the vertical planes are the only
stable entire $H$-minimal graphs}. As a corollary of our results, we
will answer affirmatively this amended conjecture in Theorem
\ref{T:instability}. We note that in Theorem \ref{T:stablevp} below
we show that the vertical planes \eqref{vp} are stable. An
$H$-minimal surface with empty characteristic locus is called
\emph{stable} if the second variation of the $H$-perimeter is
nonnegative for every compactly supported deformation, see
Definition \ref{D:minimal} below. The role of stability in the
sub-Riemannian Bernstein problem has been recently highlighted in
\cite{DGN3}, where the first three named authors have proved the
instability of the $H$-minimal entire graphs $x = y (\alpha t
+\beta)$, with $\alpha>0$, $\beta\in \R$. This result also clarified
an incorrect belief by several experts in the field, namely that
such surfaces should constitute a counterexample to the above
formulated amended form of the Bernstein conjecture.

We are ready to give a summary of our results. In $\HH$ we single
out a large class of $H$-minimal surfaces, which we call graphical
strips, see Definition \ref{D:gs} below, and which, after possibly a
left-translation and rotation about the $t$-axis, can be represented
in one of the two forms \eqref{ceI}, or \eqref{ceII}. If for the
function $G$ in these definitions we have $G'>0$ on some
sub-interval, we call the relative surface a strict graphical strip.
In Theorem \ref{T:counterexI} we show that graphical strips are
$H$-minimal, and have empty characteristic locus. Our first main
result shows that every strict graphical strip is unstable, in the
sense that there exist local deformations of the surface which
strictly increase the horizontal perimeter, see Theorem
\ref{T:perminI}. Our second main result, Theorem \ref{T:reductionI},
shows that, modulo left-translations and rotations about the group
center, every $H$-minimal entire graph in $\HH$, with empty
characteristic locus, and which is not itself a vertical plane,
contains a strict graphical strip. Combining this result with
Theorem \ref{T:perminI}, we prove that the only stable $H$-minimal
entire graphs in $\HH$, with empty characteristic locus, are the
vertical planes, see Theorem \ref{T:instability}.

To state our main theorems we begin with a definition which plays a
central role in our work.

\begin{dfn}\label{D:gs}
We say that a $C^2$ surface $\mS\subset \HH$ is a \emph{graphical
strip} if there exist an interval $I\subset \R$, and $G \in C^2(I)$,
with $G'\geq 0$ on $I$, such that, after possibly a left-translation
and a rotation about the $t$-axis, then either
\begin{equation}\label{ceI}
\mathcal S\ =\ \{(x,y,t)\in  \HH \mid  (y,t) \in \R \times I , x = y
G(t)\}\ ,
\end{equation}
or
\begin{equation}\label{ceII}
\mathcal S\ =\ \{(x,y,t)\in  \HH \mid  (x,t) \in \R \times I , y = -
x G(t)\}\ .
\end{equation}
If there exists $J\subset I$ such that $G'>0$ on $J$, then we call
$\mS$ a \emph{strict graphical strip}.
\end{dfn}

\begin{rmrk}\label{R:vp}
We stress the importance of the assumed strict positivity of $G'$ in
the definition of strict graphical strip, as opposed to the weaker
requirement $G'\geq 0$. This positivity will play a crucial role in
the proof of Theorem \ref{T:perminI}. We also mention explicitly
that while a vertical plane such as \eqref{vp} is a graphical strip,
it is not a strict graphical strip. If for instance $a\neq 0$, then
we can re-write \eqref{vp} as $x = \alpha y + \beta$, with $\alpha =
- b/a$, and $\beta = \gamma/a$. Assuming $\beta = 0$, which can be
always achieved by a left-translation (see \eqref{ltplaneI} below),
we would have $G(t) \equiv \alpha$, and therefore $G'\equiv 0$,
against the assumption in Definition \ref{D:gs}.
\end{rmrk}

We now state a first theorem which, besides having an interest in
its own right, also serves to motivate our main results.

\begin{thrm}\label{T:counterexI} Let $\mS\subset \HH$ be a graphical strip, then $\mS$ is $H$-minimal and it has empty characteristic locus.
When $I = \R$, then every surface such as \eqref{ceI} or
\eqref{ceII} is a global intrinsic $H$-minimal graph.
\end{thrm}

 We stress that in the Heisenberg group the
left-translations \eqref{Hn} are affine transformations, thereby
they preserve planes and lines. For instance, if $P$ denotes the
plane $ax + by + ct = \gamma$, then denoting by $P' = g_0 \circ P$,
where $g_0= (x_0,y_0,t_0)$, one easily sees that $P'$ is given by
\begin{equation}\label{ltplaneI}
\left(a + \frac{cy_0}{2}\right) x + \left(b - \frac{cx_0}{2}\right)
y + c t\ =\ \gamma + a x_0 + b y_0 + ct_0\ .
\end{equation}

Notice that a vertical plane ($c=0$) is mapped into a vertical plane
by a left-translation. More in general, the left-translations
preserve the property of a surface of having empty characteristic
locus. Furthermore, they preserve the $H$-mean curvature, and
therefore the $H$-minimality, the $H$-perimeter, and the property of
a surface of being stable. We also notice that rotations about the
$t$-axis (the group center), also have the same properties.

The notion of intrinsic graph in the second part of Theorem
\ref{T:counterexI} is that introduced in \cite{FSS4}, but see also
\cite{FSS2}. The proof of Theorem \ref{T:counterexI} shows that if
$\mS$ is of type \eqref{ceI} with $I= \R$, then it is a global
$X_1$-graph, whereas if it is of type \eqref{ceII}, then $\mS$ is a
global $X_2$-graph. We recall that $\mS$ is said an intrinsic
$X_1$-graph if there exist an open set $\Om \subset \R^2_{u,v}$, and
a $C^2$ function $\phi:\Om \to \R$, such that on $\Om$ we can
parameterize $\mS$ as follows $(x,y,t) = (0,u,v) \circ \phi(u,v) e_1
= (0,u,v) \circ (\phi(u,v),0,0)$. This means that
\begin{equation}\label{X1graph}
\mS\ =\ \left\{(x,y,t)\in \HH\mid (u,v)\in \Om, (x,y,t) =
\left(\phi(u,v),u,v - \frac{u}{2} \phi(u,v)\right)\right\}\ .
\end{equation}

When we can take $\Om =\R^2_{(u,v)}$, then $\mS$ is called a global
$X_1$-graph. Similarly, if the points of $\mS$ can be described by
$(x,y,t) = (u,0,v) \circ \phi(u,v) e_2 = (u,0,v) \circ
(0,\phi(u,v),0)$, i.e., if
\begin{equation}\label{X2graph}
\mS\ =\ \left\{(x,y,t)\in \HH\mid (u,v)\in \Om, (x,y,t) =
\left(u,\phi(u,v),v + \frac{u}{2} \phi(u,v)\right)\right\}\ ,
\end{equation}
then $\mS$ is said an intrinsic $X_2$-graph (global, if $\Om=
\R^2_{(u,v)}$). Clearly, the vertical planes \eqref{vp} are global
intrinsic graphs. If, for instance, $a\neq 0$, then $\tilde P_0$ can
be parameterized as in \eqref{X1graph}, with
\[
\phi(u,v)\ =\ -\ \frac{b}{a} u + \frac{\gamma}{a}\ .
\]

 Before proceeding, we pause to give some
examples which illustrate the situation of Theorem
\ref{T:counterexI}.

\textbf{Example 1:} The choice $G(t) = \tan \tanh(t)$ makes
\eqref{ce} a special case of Theorem \ref{T:counterexI}. We stress
that in this example $G'(t)>0$ for every $t\in \R$, and therefore
the surface $\mS$ is a strict graphical strip with $I = \R$.
According to Theorem \ref{T:counterexI}, we conclude that $\mS$ is
also a $C^\omega$ global intrinsic $X_1$-graph.

\textbf{Example 2:} The class of strict graphical strips is not
contained in that of global intrinsic graphs. Consider, for
instance, the function $G(t) = \cot(- t + \frac{\pi}{2})$, with $-
\frac{\pi}{2} < t < \frac{\pi}{2}$. We note that $G'(t)>0$, and
therefore the corresponding surface in \eqref{ceI},
\[ x\ =\ y\ \cot(- t + \frac{\pi}{2})\ ,
\]
is a strict graphical strip, with $I  = (-
\frac{\pi}{2},\frac{\pi}{2})$. This surface is the helicoid
\[
x\ =\ r\ \cos \theta\ ,\ y\ =\ r\ \sin \theta\ , \ t\ =  -  \theta +
\frac{\pi}{2} \ ,\ 0\leq r< \infty\ , \ 0<\theta<\pi \ .
\]

If we set $D = \R\times I$, then the map $\Phi(y,t) = (y,t +
\frac{y^2}{2}G(t))$ is not a diffeomorphism of $D$ onto the whole
$(u,v)$-plane. However, the arguments in the proof of Theorem
\ref{T:counterexI} show that it is a diffeomorphism onto the
connected open subset $\Om\subset \R^2_{u,v}$ defined by $\Om = \R^2
\setminus (L^+ \cup L^-)$, where $L^+ = \{(0,v)\in \R^2\mid v\geq
\frac{\pi}{2}\}$, $L^- = \{(0,v)\in \R^2\mid v\leq -
\frac{\pi}{2}\}$. As a consequence, the helicoid is \emph{not} a
global intrinsic graph, although it is a $C^\omega$ intrinsic
$X_1$-graph on the domain $\Om$. If we denote by $\Psi(u,v) =
(\Psi_1(u,v),\Psi_2(u,v)) = \Phi^{-1}(u,v)$, then the function
$\phi(u,v)$ in \eqref{X1graph} is given by
\[
\phi(u,v)\ =\ u\ G(\Psi_2(u,v))\ ,\quad\quad\quad (u,v)\in \Om\ .
\]

\textbf{Example 3:} The choice $G(t) = \alpha t + \beta$, with
$\alpha>0$, $\beta\in \R$, gives the strict graphical strips $x =
y(\alpha t + \beta)$ studied in \cite{DGN3}, where it was proved
that these surfaces are unstable. We note that such surfaces are
$C^\omega$ global intrinsic $H$-minimal graphs . We also observe
that in this example it is possible to compute explicitly the
function $\phi(u,v)$ in \eqref{X1graph} which describes $\mS$ as a
global $X_1$-graph. One finds
\begin{equation}\label{phi1}
\phi(u,v)\ =\ \frac{2 u (\alpha v + \beta)}{2 + \alpha u^2}\ .
\end{equation}

\textbf{Example 4:} The surface described by \[ \mS\ =\ \{(x,y,t)\in
\HH \mid (y,t)\in \R\times (0,\infty) , x = y t^2\}\ , \] is a
strict graphical strip with $I = (0,\infty)$, but like the surface
in example 2, it is not a global intrinsic graph.

 We are now ready to state the first main result of this paper.
In what follows we indicate with $\nuX = \bN^H/|\bN^H|$ the
horizontal Gauss map of $\mS$, see section \ref{S:prelim}, and by
$\mathcal V^H_{II}(\mS; \mathcal X)$ the second variation of the
$H$-perimeter with respect to a deformation of $\mS$ in the
direction of the vector field $\mathcal X$, see Definitions
\ref{D:variations} and \ref{D:minimal} below. An $H$-minimal surface
$\mS$ is called \emph{stable} if $\sv \geq 0$ for every compactly
supported $\mathcal X = a X_1 + b X_2 + k T$. Otherwise, it is
called unstable. We note that, since thanks to Theorem
\ref{T:counterexI} every graphical strip has empty characteristic
locus, the horizontal Gauss map $\nuX$ of such a surface is globally
defined.

\begin{thrm}\label{T:perminI}
Let $\mS$ be a strict graphical strip, then $\mS$ is unstable. In
fact, there exists a continuum of $h\in C^2_0(\mS)$, for which
$\mathcal V^H_{II}(\mS; h \nuX) <  0$.
\end{thrm}

As a consequence of Theorem \ref{T:perminI}, the $H$-minimal
surfaces corresponding to Examples 1, 2, 3 and 4, are all unstable,
i.e., they are not local minimizers of the $H$-perimeter. We
emphasize that the class of strict graphical strips is very wide.
For instance, as we show in Theorem \ref{T:reductionI} below, every
$H$-minimal entire graph in $\HH$, with empty characteristic locus,
and which is not itself a vertical plane, contains a strict
graphical strip. The main ingredients in the proof of Theorem
\ref{T:perminI} are the second variation formulas for the
$H$-perimeter, see Theorem \ref{T:nudef} below, and the explicit
construction of deformations of the surface along which the
$H$-perimeter decreases strictly. This part is delicate and it has
been influenced by the recent paper \cite{DGN3}. In connection with
Theorem \ref{T:perminI} we mention that, after the present work was
completed, we have received the interesting paper \cite{BSV} in
which the authors, using the construction in \cite{DGN3}, in
combination with other tools in part also independently developed in
\cite{GS}, establish the instability of global intrinsic graphs,
thus answering affirmatively the above formulated Bernstein type
conjecture in this setting. As we have seen, Theorem \ref{T:perminI}
includes surfaces, such as for instance the helicoid in example 2,
or that in example 4, which are not global intrinsic graphs.

Here is our second main result.

\begin{thrm}\label{T:reductionI} Let
$\mS\subset \HH$ be an $H$-minimal entire graph, with empty
characteristic locus, and that is not itself a vertical plane such
as \eqref{vp}, then there exists a strict graphical strip $\mS_0
\subset \mS$.
\end{thrm}

We mention explicitly that the above statement should be interpreted
in the sense that, after possibly composing with a suitable rotation
about the $t$-axis and a left-translation, the transformed surface
contains a portion of the type \eqref{ceI}, or \eqref{ceII}. The
proof of Theorem \ref{T:reductionI} is based on some of the main
results in \cite{GP}. For clarity of exposition we present in
Section \ref{S:rs} a detailed account of the main reduction steps.
We mention that, alternatively, one could use the results
independently obtained in \cite{CHMY} and \cite{CH}.

Finally, by combining Theorems \ref{T:perminI} and
\ref{T:reductionI}, we answer affirmatively the above formulated
Bernstein conjecture.

\begin{thrm}[\textbf{of Bernstein type}]\label{T:instability}
In $\HH$ the only stable $H$-minimal entire graphs, with empty
characteristic locus, are the vertical planes \eqref{vp}.
\end{thrm}

Concerning the higher-dimensional case of Theorem
\ref{T:instability} we mention that, using the construction in
\cite{BDG}, we obtain the following negative result.

\begin{thrm}\label{T:negative}
In $\Hn$, with $n\geq 5$, there exist $C^\omega$ stable $H$-minimal
graphs, with empty characteristic locus, which are not vertical
hyperplanes.
\end{thrm}

What happens when $n=2,3,4$ is, presently, terra incognita. We will
now briefly describe the organization of the paper. In section
\ref{S:prelim} we collect some basic facts about the Heisenberg
group, and introduce the main geometric set-up. In section
\ref{S:1&2var} we collect some results from sub-Riemannian geometric
measure theory, specifically, the first and second variation
formulas for the horizontal perimeter. Sections \ref{S:proof} and
\ref{S:rs} are the central parts of the paper. After proving Theorem
\ref{T:counterexI}, the remainder of section \ref{S:proof} is
devoted to proving Theorem \ref{T:perminI}. In section \ref{S:rs} we
prove Theorems \ref{T:reductionI} and \ref{T:instability}. Finally,
in section \ref{S:bdg} we prove Theorem \ref{T:negative}.

\medskip

\noindent \textbf{Acknowledgment:} This paper was presented by the
second named author in the UCI-UCSD joint Differential Geometry
Seminar at the University of California, San Diego, La Jolla, in
February 2006. He would like to thank Lei Ni for stimulating
discussions and for his gracious hospitality during his visit, and
Salah Baouendi, Dima Khavinson, Peter Li and Linda Rothschild for
their interest in the problem.

\vskip 0.6in


\section{\textbf{Preliminaries}}\label{S:prelim}

\vskip 0.2in

In this section we collect some definitions and known results which
will be needed in the paper. We recall that the Heisenberg group
$\Hn$ is the graded, nilpotent Lie group of step $r=2$ whose
underlying manifold is $\mathbb C^{n}\times \R \cong \R^{2n+1}$,
with non-Abelian left-translation
\[
L_{(z,t)} (z',t')\ =\ (z,t) \circ (z',t')\ =\ \left(z + z', t + t'
-\frac{1}{2} Im (z \cdot \overline{z'})\right)\ ,
\]
and non-isotropic dilations
\begin{equation}\label{Hn3}
\delta_\la (z,t)\ =\ (\la z, \la^2 t)\ ,\quad\quad\quad \la > 0\ .
\end{equation}

These dilations provide a natural scaling associated with the
grading of the Heisenberg algebra $\mathfrak h_n = V_1\oplus V_2$,
where $V_1 = \mathbb C^n \times \{0\}$, $V_2 = \{0\}\times \R$. The
homogeneous dimension associated with \eqref{Hn3} is $Q = 2n + 2$.
We recall that identifying $\mathfrak h_n$ with $\mathbb R^{2n+1}$,
by identifying $z = x + i y \in \mathbb C^n$ with $(x,y)\in
\R^{2n}$, we have for the bracket of $\xi = (x,y,t)$, $\xi'
=(x',y',t') \in\mathfrak h_n$
\[
[\xi,\xi']\ =\ (0,0, x\cdot y'- x'\cdot y)\ .
\]

Here, and throughout the paper, we will use $v \cdot w$ to denote
the standard Euclidean inner product of two vectors $v$ and $w$ in
$\Rn$. It is then clear that $[V_1,V_1] = V_2$, and that $V_2$ is
the group center. Via the Caley map, $\Hn$ can be identified with
the boundary of the Siegel upper half-space $\mathcal U^n = \{z\in
\mathbb C^{n+1}\mid Im\ z_{n+1}
> 2 \sum_{j=1}^n |z_j|^2\}$, see Ch.12 in \cite{S}.
In the  real coordinates $g =(x,y,t)\in \R^{2n+1}$ the non-Abelian
group law of $\Hn$ is given by
\begin{equation}\label{Hn}
g \circ g'\ =\ (x,y,t) \circ (x',y',t')\ =\ (x + x', y+ y', t+t' +
\frac{1}{2} (x\cdot y' - x' \cdot y))\ . \end{equation}

Let $(L_g)_*$ be the differential of the left-translation $L_g(g')
= g \circ g'$. A simple computation shows that
\begin{align}\label{Hn2}
& (L_g)_*\left(\frac{\partial}{\partial x_i}\right)\
\overset{def}{=}\ X_i\ =\ \frac{\partial}{\partial x_i}\ -\
\frac{y_i}{2}\ \frac{\partial}{\partial t}\ , \quad\quad i =
1,...,n\ ,
\\
& (L_g)_*\left(\frac{\partial}{\partial y_i}\right)\
\overset{def}{=}\ X_{n+i}\ =\ \frac{\partial}{\partial y_i}\ +\
\frac{x_i}{2}\ \frac{\partial}{\partial t}\ , \quad\quad i =
1,...,n\ ,
\notag\\
& (L_g)_*\left(\frac{\partial}{\partial t}\right)\
\overset{def}{=}\ T\ =\ \frac{\partial}{\partial t} \notag
\end{align}

We note that the only non-trivial commutator is
\[
[X_i,X_{n+j}]\ =\ \delta_{ij}\ T\ , \quad\quad\quad i, j =
1,...,n\ ,
\]
therefore the vector fields $\{X_1,...,X_{2n}\}$ generate the Lie
algebra $\mathfrak h_n$.

Henceforth, $\Hn$ will be endowed with a left-invariant inner
product $<\cdot,\cdot>$, with respect to which
$\{X_1,...,X_{2n},T\}$ constitute an orthonormal basis. With the
exception of the Euclidean inner product, which, as we have said
above, we denote $v \cdot w$, no other inner product will be used in
this paper, so when we write $<\cdot,\cdot>$ there will be no danger
of confusion. The horizontal bundle is $H\Hn = \cup_{g\in \Hn} H_g
\Hn$, where \[ H_g\Hn\ =\ span\{X_1(g),...,X_{2n}(g)\}\ . \]

We denote by $\nabla_X Y$ the Levi-Civita connection with respect to
$<\cdot,\cdot>$. Projecting such connection onto the horizontal
subbundle $H\Hn \subset T\Hn$, we obtain a connection $\nabla^H_X Y$
on $H\Hn$, which we call the \emph{horizontal Levi-Civita
connection}. This idea goes back to that of E. Cartan's
\emph{non-holonomic connection}, see \cite{C}. For any $X\in
\Gamma(T\Hn)$, $Y\in \Gamma(H\Hn)$ we let
\begin{equation}\label{hlc}
\nabla^H_X Y\ =\ \sum_{i=1}^{2n} <\nabla_X Y,X_i> X_i\ ,
\end{equation}
and one can easily verify that $\nabla^H_X Y$ is metric preserving
and torsion free, in the sense that if we define the horizontal
torsion of $\mS$ as
 \[
 T^H(X,Y)\ =\ \nabla^H_X Y\ -\ \nabla^H_Y X\ -\ [X,Y]^H\ ,
 \]
where $[X,Y]^H = \sum_{i=1}^{2n} <[X,Y],X_i>X_i$, then $T^H(X,Y)=
0$. If $f\in C^1(\Hn)$, we let
\[
\nabla^H f\ =\ \sum_{i=1}^{2n} <\nabla f,X_i> X_i\ ,
\]
where we have denoted by $\nabla f$ the Riemannian gradient of
$f$.

Given an oriented hypersurface $\mS \subset \Hn$, we denote by $\bN$
the Riemannian non-unit normal to $\mS$. Throughout the paper, we
will indicate with
\[
\Sigma(\mS)\ =\ \{g\in \mS \mid T_g\mS = H_g \Hn\}\ , \]
characteristic locus of $\mS$. Notice that $g\in \Sigma(\mS)$ is
equivalent to having $<\bN,X_i> = 0$, $i=1,...,2n$, at $g$. We
recall that it was proved in \cite{B} that $H^{Q-1}(\Sigma(\mS))=
0$, where $H^s$ indicates the $s$-dimensional Hausdorff measure
constructed with the Carnot-Carath\'eodory distance associated with
the subbundle $H\Hn$. If $\mS$ is of class $C^2$, then we define the
non-unit horizontal normal of $\mS$ by
\[
\bN^H\ =\ \sum_{i=1}^{2n} <\bN,X_i>X_i\ .
\]

It is clear that at a given $g\in \mS$ one has $\bN^H \not= 0$ if
and only if $g\not\in \Sigma(\mS)$. If $g\not\in \Sigma(\mS)$, then
the horizontal tangent space to $\mS$ in $g$ is defined by
\begin{equation}\label{hts}
HT_g \mS\ =\ \{\boldsymbol v\in H_g\Hn \mid <\boldsymbol v,\bN^H> =
0\}\ .
\end{equation}

It is easy to recognize that $HT_g \mS = T_g \mS \cap H_g\Hn$. The
horizontal tangent bundle of $\mS$ is $HT \mS = \underset{g\in
\mS}{\cup}T^H_g \mS$. At every non-characteristic point of $\mS$ we
define the horizontal Gauss map by letting
\[
\nuX\ =\ \frac{\bN^H}{|\bN^H|}\ .
\]

Since
\[
H_g\Hn\ =\ HT_g\mS\ \oplus\ span\{\nuX_g\}\ ,
\]
we have $dim(HT_g\mS) = 2n -1$. For instance, in $\HH$ we simply
have $HT\mS = span\{\boldsymbol e_1\}$, where
\[
\boldsymbol e_1\ =\ \nup\ =\ <\nuX,X_2> X_1\ -\ <\nuX,X_1> X_2\ .
\]

\medskip

\begin{dfn}\label{D:horconS}
Let $\mS\subset \Hn$ be a $C^k$ hypersurface, $k\geq 2$, with
$\Sigma(\mS) = \varnothing$, then we define the \emph{horizontal
connection} on $\mS$ as follows. Let $\nabla^H$ denote the
horizontal Levi-Civita connection introduced above. For every
$X,Y\in C^1(\mS;HT\mS)$ we let
\[
\delXY\ =\ \nabla^H_{\oX} \oY\ -\ <\nabla^H_{\oX} \oY,\nuX> \nuX\
,
\]
where $\oX, \oY \in C^1(\Hn;H\Hn)$ are such that $\oX = X$, $\oY =
Y$ on $\mS$.
\end{dfn}

\medskip

One can check that Definition \ref{D:horconS} is well-posed, i.e.,
it is independent of the extensions $\oX, \oY$ of the vector
fields $X, Y$. For every $X,Y\in C^1(\mS;HT\mS)$ one has
\begin{equation}\label{delXY}
\delXY\ -\ \nabla^{H,\mS}_Y X\ =\ [X,Y]^H\ -\ <[X,Y]^H,\nuX>\nuX\
,
\end{equation}
in other words $\delXY - \delYX$ equals the projection of
$[X,Y]^H$ onto $HT\mS$.

It is clear from \eqref{delXY} that the horizontal connection
$\nabla^{H,\mS}$ on $\mS$ is not necessarily torsion free. This
depends on the fact that it is not true in general that, if
$X,Y\in C^1(S;HT\mS)$, then $[X,Y]^H \in C^1(\mS;HT\mS)$. However,
since in the first Heisenberg group $\HH$ this is trivially true,
in this setting $\nabla^{H,\mS}$ is torsion free, and therefore it
has the properties of a Levi-Civita connection.

Given a function $f\in C^1(\mS)$ we will denote
\begin{equation}\label{hg}
\nabla^{H,\mS} f\  =\ \nabh \overline f\ -\ <\nabh \overline
f,\nuX>\nuX\ ,
\end{equation}
where $\overline f\in C^1(\Hn)$ denotes any extension of $f$.
Henceforth, we will let
\[
\di f \ =\ <\nabla^{H,\mS} f, X_i>\ =\ X_i \overline f\ -\ <\nabh
\overline f,\nuX> <\nuX,X_i>\ .
\]

\medskip

\begin{dfn}\label{D:sff}
Let $\mS\subset \Hn$ be a $C^k$ hypersurface, $k\geq 2$, with
$\Sigma(\mS) = \varnothing$, then for every $X,Y\in C^1(\mS;HT\mS)$
we define a tensor field of type $(0,2)$ on $\mS$, as follows
\begin{equation}\label{sff} II^{H,\mS}(X,Y)\ =\ <\nabla^H_X Y,\nuX>
\nuX\ .
\end{equation}
We call $II^{H,\mS}(\cdot,\cdot)$ the \emph{horizontal second
fundamental form} of $\mS$. We also define $\mathcal A^{H,\mS} :
HT \mS \to HT \mS$ by letting for every $g\in \mS$ and
$\boldsymbol u, \boldsymbol v \in HT_{g}\mS$
\begin{equation}\label{shape}
<\mathcal A^{H,\mS} \boldsymbol u,\boldsymbol v>\ =\ -\
<II^{H,\mS}(\boldsymbol u,\boldsymbol v),\nuX>\ =\ -\ <\nabla_X^H
Y,\nuX>\ ,
\end{equation}
where $X, Y \in C^1(\mS,HT \mS)$ are such that $X_g = \boldsymbol
u$, $Y_g = \boldsymbol v$. We call the linear map $\mathcal
A^{H,\mS} : HT_{g}\mS \to HT_{g}\mS$ the \emph{horizontal shape
operator} If $\boldsymbol e_1,...,\boldsymbol e_{2n-1}$ denotes a
local orthonormal frame for $HT\mS$, then the matrix of the
horizontal shape operator with respect to the basis $\boldsymbol
e_1,...,\boldsymbol e_{2n-1}$ is given by the $(2n-1)\times(2n-1)$
matrix \ $- \big[<\nabla_{\boldsymbol e_i}^H \boldsymbol
e_j,\nuX>\big]_{i,j=1,...,2n-1}$.

\end{dfn}

\medskip

Definitions \ref{D:horconS} and \ref{D:sff} are taken from
\cite{DGN1}. A different notion of the second fundamental form has
been explored by the last named author and R. Hladky [HP] using a
generalization of the Webster-Tanaka connection on a wide class of
sub-Riemannian manifolds.  This second fundamental form and the
(un-symmetrized) operator $\mathcal A^{H,\mS}$ are used to analyze
the minimal and constant mean curvature surfaces in this setting. We
emphasize that, when restricted to the case of the Carnot groups,
these two formulations are equivalent.

We call \emph{horizontal principal curvatures} the real eigenvalues
$\kappa_1,...,\kappa_{2n-1}$ of the symmetrized operator $\mathcal
A^{H,\mS}_{sym} = \frac{1}{2}(\mathcal A^{H,\mS} + (\mathcal
A^{H,\mS})^t)$. The \emph{horizontal mean curvature} of $\mS$ is
defined as follows
\begin{equation}\label{hmc}
\mathcal H\ =\ \kappa_1\ +\ ...\ +\ \kappa_{2n-1}\ .
\end{equation}

When the hypersurface $\mS$ has non empty characteristic locus,
then at every $g_0\in \Sigma(\mS)$ we define
\[
\mathcal H(g_0)\ =\ \underset{g\to g_0, g\in \mathcal S\setminus
\Sigma}{\lim}\ \mathcal H(g)\ ,
\]
provided that such limit exists, finite or infinite. We do not
define the $H$-mean curvature at those points $g_0\in \Sigma$ at
which the above limit does not exist.

\medskip

\begin{dfn}\label{D:minimal}
A $C^2$ surface $\mS\subset \HH$ is called $H$-minimal if
$\mathcal H \equiv 0$ as a continuous function on $\mS$.
\end{dfn}

\medskip

We will need the following result, which will prove useful for
computing the $H$-mean curvature, see Proposition 9.8 in
\cite{DGN1}.

\medskip

\begin{prop}\label{P:equalMC}
The $H$-mean curvature in definition \eqref{hmc} coincides with
the function defined by the equation
\[
\mathcal H\ =\ \sum_{i=1}^{2n}\ \di\ <\nuX,X_i>\ .
\]
\end{prop}

\medskip

Given a vector field $X\in C^1(\Hn;H\Hn)$, we denote by
\[
div^H X\ =\ \sum_{i=1}^{2n} X_i <X,X_i>\ ,
\]
the horizontal divergence of $X$.

\vskip 0.6in

\section{\textbf{First and second variation of the
$H$-perimeter}}\label{S:1&2var}

\vskip 0.2in

In this section we introduce the relevant notions of stationary and
stable surface which enter in the statement of Theorem
\ref{T:instability}, and we recall the first and second variation
formulas from \cite{DGN1} and \cite{DGN3} which will be used in its
proof.

Given an oriented $C^2$ surface $S\subset \HH$, with Riemannian
(non-unit) normal $\bN$, we introduce the quantities
\begin{equation}\label{pq}
p\ =\ <\bN, X_1>\ ,\ \quad q\ =\ <\bN , X_2>\ ,\quad\quad \omega\ =\
<\bN,T>\ ,\ \quad W\ =\ \sqrt{p^2 + q^2}\ ,
\end{equation}
and, at every point where $W\neq 0$, we set
\begin{equation}\label{bars}
\pb\ =\ \frac{p}{W}\ ,\quad\quad \qb\ =\ \frac{q}{W}\ ,\quad\quad
\ob\ =\ \frac{\omega}{W}\ .
\end{equation}

Notice that
\begin{equation}\label{normals}
\bN^H\ =\ p X_1 + q X_2\ ,\quad\quad \nuX\ =\ \pb X_1 + \qb X_2\
,\quad\quad <\bN^H,\bN>\ =\ W^2\ .
\end{equation}

>From \eqref{normals} we easily recognize that
\begin{equation}\label{cosine}
\cos(\bN^H\angle \bN)\ =\ \frac{W}{|\bN|}\ .
\end{equation}

In the classical theory of minimal surfaces, the concept of area or
perimeter occupies a central position, see \cite{DG1}, \cite{DG2},
\cite{G}, \cite{MM}, \cite{Si}, \cite{CM}. In sub-Riemannian
geometry there exists an appropriate variational notion of
perimeter. Given an open set $\mathcal U \subset \HH$ we denote
$\mathcal F(\mathcal U) = \{\zeta \in C^1_0(\mathcal U;H\HH)\mid
||\zeta||_{L^\infty(\mathcal U)}\leq 1\}$. A function $f\in
L^1(\mathcal U)$ is said to belong to $BV_H(\mathcal U)$ (the space
of functions with finite horizontal bounded variation), if
\[
Var_H(f;\mathcal U)\ =\ \underset{\zeta\in \mathcal F(\mathcal
U)}{\sup}\ \int_\mathcal U f\ div^H \zeta\ dg\ < \ \infty\ .
\]

This space becomes a Banach space with the norm $||f||_{BV_H(\Om)}
= ||f||_{L^1(\Om)} + Var_H(u;\Om)$. Given a measurable set
$\mathcal E\subset \HH$, the $H$-\emph{perimeter} of $\mathcal E$
with respect to the open set $\mathcal U\subset \HH$ is defined as
follows, see for instance \cite{CDG}, and \cite{GN},
\[
P_H(\mathcal E;\mathcal U)\ =\ Var_H(\chi_{\mathcal E};\mathcal
U)\ .
\]

When $\mathcal E$ is a $C^1$ domain, one can recognize that
\begin{equation}\label{per}
P_H(\mathcal E;\mathcal U)\ =\ \int_{\p \mathcal E\cap \mathcal U}
\cos(\bN^H \angle \bN)\ d\sigma\ ,
\end{equation}
where $d\sigma$ indicates the standard surface measure on $\p
\mathcal E$. We will denote by $d\sigma_H$ the $H$-perimeter measure
concentrated on $\mS=\p \mathcal E$. According to \eqref{cosine},
\eqref{per}, we have for any Borel subset $E\subset \mS$ such that
$P_H(E)<\infty$,
\begin{equation}\label{pm}
\sigma_H(E)\  =\ \int_E \frac{W}{|\bN|}\ d\sigma\ .
\end{equation}

Two important properties of the $H$-perimeter are its invariance
with respect to the dilations \eqref{Hn3} and the left-translations
\eqref{Hn}. The former, is expressed by the equation
\[
\sigma_H(\delta_\lambda(\mS))\ =\ \lambda^{Q-1}\ \sigma_H(\mS)\ .
\]

In keeping up with the notation of \cite{DGN3} it will be
convenient to indicate with $Y\zeta$ and $Z\zeta$ the respective
actions of the vector fields $\nuX$ and $\nup$ on a function
$\zeta\in C^1(\mathcal S)$, thus
\begin{equation}\label{Y}
Y\zeta\ \overset{def}{=}\ \pb\ X_1 \zeta\ +\ \qb\ X_2 \zeta\
,\quad\quad\quad Z\zeta\ \overset{def}{=}\ \qb\ X_1\zeta\ -\ \pb\
X_2\zeta\ .
\end{equation}

The frame $\{Z,Y,T\}$ is orthonormal. It is worth observing that,
since the metric tensor $\{g_{ij}\}$ with respect to the inner
product $<\cdot,\cdot>$ has the property $det\{g_{ij}\} = 1$, then
the (Riemannian) divergence in $\HH$ of these vector fields is
given by
\begin{equation}\label{divY}
div\ Y\ =\ X_1 \pb\ +\ X_2 \qb\ =\ \mathcal H\ ,\quad\quad\quad
div\ Z\ =\ X_1 \qb\ -\ X_2 \pb\ ,
\end{equation}
where the first equality is justified by Proposition
\ref{P:equalMC} and by the fact that $|\nuX| = 1$. Using Cramer's
rule one easily obtains from \eqref{Y}
\begin{equation}\label{Xs}
X_1 \zeta\ =\ \qb\ Z\zeta\ +\ \pb\ Y\zeta\ ,\quad\quad\quad
X_2\zeta\ =\ \qb\ Y\zeta\ -\ \pb\ Z\zeta\ .
\end{equation}

One also has
\begin{equation}\label{deltas}
\duno \zeta\ =\ \qb\ Z\zeta\ ,\quad\quad\quad\quad \ddue \zeta\ =\
-\ \pb\ Z\zeta\ ,
\end{equation}
so that
\begin{equation}\label{delZ}
|\nabla^{H,\mS} \zeta|^2\ =\ (Z\zeta)^2\ .
\end{equation}

We notice that
\begin{equation}\label{mc}
\qb Z\pb\ -\ \pb Z\qb\ =\ \mathcal H\ . \end{equation}

This can be easily recognized using the first equation in
\eqref{divY}, and \eqref{Xs}, as follows
\[ \mathcal H\ =\ X_1
\pb\ +\ X_2 \qb\ =\ \qb Z\pb - \pb Z\qb + \pb Y\pb + \qb Y\qb\ =\
\qb Z\pb - \pb Z\qb\ ,
\]
where we have used the fact that $0 = \frac{1}{2} Y(\pb^2 + \qb^2)
= \pb Y\pb + \qb Y\qb$.

\medskip

\begin{dfn}\label{D:variations}
Let $\mathcal S \subset \HH$ be an oriented $C^2$ surface, with
$\Sigma(\mS)=\varnothing$. Consider the family of vector fields
$\mathcal X = a X_1 + b X_2 + k T$, with $a, b , k\in C^2_0(\mathcal
S)$, and the family of surfaces $\mS^\la$, where for small $\lambda
\in \mathbb R$ we have let
\begin{equation}\label{def}
\mathcal S^\lambda = J_\lambda(\mathcal S)\ =\ \mS\ +\ \la
\mathcal X\ .
\end{equation}
We define the \emph{first variation} of the $H$-perimeter with
respect to the deformation \eqref{def} as
\[
\fv\ =\ \frac{d}{d\lambda}~ P_H(\mathcal S^\lambda)\Bigl|_{\lambda
= 0}\ .
\]
We say that $\mS$ is \emph{stationary} if $\fv = 0$, for every
$\mathcal X$.
\end{dfn}

\medskip

Classical minimal surfaces are stationary points of the perimeter
(the area functional for graphs). It is natural to ask what is the
connection between the notion of $H$-minimal surface and that of
$H$-perimeter. To answer this question we recall the following
results from \cite{DGN1}.

\medskip

 \begin{thrm}\label{T:variations}
Let $\mS\subset \HH$ be an oriented $C^2$ surface, with $\Sigma(\mS)
= \varnothing$, then
\begin{equation}\label{fvH}
\mathcal V^H_I(\mS;\mathcal X)\ =\
 \int_{\mathcal S}
\mathcal H\ \frac{\cos(\mathcal X \angle \bN)}{\cos(\nuX \angle
\bN)}\ |\mathcal X|\ d\sigma_H\ .
\end{equation}
In particular, $\mathcal S$ is stationary if and only if it is
$H$-minimal.
\end{thrm}

\medskip

Versions of Theorem \ref{T:variations} have also been obtained
independently by other people. An approach based on motion by
$H$-mean curvature can be found in \cite{BC}. When $\mathcal X = a
\nuX + k T$, then a proof based on CR-geometry can be found in
\cite{CHMY}, and \cite{RR2}. Recently, a general first variation
formula for a wide class of sub-Riemannian spaces has been found in \cite{HP}.

\medskip

\begin{dfn}\label{D:minimal}
Given an oriented $C^2$ surface $\mS \subset \HH$, with $\Sigma(\mS)
= \varnothing$, we define the \emph{second variation} of the
$H$-perimeter with respect to the deformation \eqref{def} as
\[
\sv\ =\ \frac{d^2}{d\lambda^2}~ P_H(\mathcal
S^\lambda)\Bigl|_{\lambda = 0}\ .
\]
We say that $\mS$ is \emph{stable}, if it is stationary (i.e.,
$H$-minimal), and if
\[
\sv\ \geq \ 0\ ,\quad\quad\quad\text{for every}\quad \mathcal X  .
\]
If there exists $\mathcal X \neq 0$ such that $\sv < 0$, then we say
that $\mS$ is \emph{unstable}.
\end{dfn}

\medskip

The following second variation formula from \cite{DGN1}, see also
Theorem 3.3 in \cite{DGN3}, will play a crucial role in the proof of
Theorem \ref{T:perminI}.

\medskip

\begin{thrm}\label{T:nudef}
Let $\mS\subset \HH$ be a $C^2$ oriented surface, with empty
characteristic locus, then the second variation of the
$H$-perimeter with respect to the deformation \eqref{def}, with
$\mathcal X = h \nuX$, $h\in C^2_0(\mS)$, is given by
\begin{align}\label{2varnu}
& \mathcal V^H_{II}(\mS;h\nuX)\ =\ \int_\mS \bigg\{(Zh)^2 + h^2
\big[2 (\pb T\qb - \qb T\pb) + 2 \ob (\qb Y\pb - \pb Y\qb) +
\ob^2\big]\bigg\}\ d\sigma_H\ .
\end{align}
If instead we choose $\mathcal X = a X_1$, $a\in C^2_0(\mS)$, in
\eqref{def}, then the corresponding second variation is given by
\begin{align}\label{X1def}
\mathcal V^H_{II}(\mS;a X_1)\ & =\ \int_{\mathcal S} \bigg\{\pb^2
(Za)^2 + \pb^2\ \ob^2\ a^2
\\
& +\ \ob Z(a^2) - \pb\ \qb \left(T(a^2) - \ob
Y(a^2)\right)\bigg\}\ d\sigma_H\ . \notag
\end{align}
\end{thrm}

\medskip

Theorem \ref{T:nudef} ia special case of a general second variation
formula found in \cite{DGN1}, see also Theorem 3.3 in \cite{DGN3}.
Using such general result, in combination with some integration by
parts formulas from \cite{DGN1}, we can prove the following theorem.

\medskip

\begin{thrm}\label{T:stablevp}
Every vertical plane such as \eqref{vp} is stable.
\end{thrm}

\begin{proof}[\textbf{Proof}]
Since the notion of stability is invariant under left-translations,
and rotations about the $t$-axis, we can assume without restriction
the $\tilde P_0 = \mS$ is the plane $x=0$. Since for this surface we
have $\pb \equiv 1, \qb = \ob  \equiv 0$, the formula (3.3) in
Theorem 3.3 in \cite{DGN3} gives for $\mathcal X = a X_1 + b X_2 + k
T$,
\[
\sv\ =\ \int_\mS \left\{(Za)^2 + 2 (Tb Zk - Tk Zb) + T(ab)\right\}
d\sigma_H\ .
\]

Since $a,b, c$ are compactly supported in $\mS$, Lemma 3.7 in
\cite{DGN3} gives
\[
\int_\mS T(a b) d\sigma_H\ =\ 0\ .
\]

We also have
\begin{align*}
\int_\mS Tb Zk\ d\sigma_H\ &  =\ \int_\mS T(b Zk)\ d\sigma_H\ -\
\int_\mS b T(Zk)\ d\sigma_H
\\
& =\ -\ \int_\mS b Z(Tk)\ d\sigma_H\ +\ \int_\mS b [Z,T]k\ d\sigma_H
\end{align*}

We now using the following commutator formula from \cite{DGN1}
\begin{equation*}
[Z,T] \ =\ (\qb T \pb - \pb T \qb)\ Y\ ,
\end{equation*}
valid on any non-characteristic surface. In the present case, such
formula gives $[Z,T]= 0$, and we thus conclude
\begin{align*}
\int_\mS Tb Zk\ d\sigma_H\ & =\ -\ \int_\mS b Z(Tk)\ d\sigma_H\ =\
-\ \int_\mS Z(b Tk)\ d\sigma_H\ +\ \int_\mS Tk Zb\ d\sigma_H\ .
\end{align*}

Since Lemma 3.6 in \cite{DGN3} gives for every $\zeta\in
C_0^1(\mS)$,
\[
\int_{\mathcal S} Z\zeta\ d\sigma_H\ =\ -\ \int_{\mathcal S} \zeta\
\ob\  d\sigma_H\ =\ 0\ ,
\]
we finally obtain
\[
\sv\ =\ \int_\mS (Za)^2  d\sigma_H\ \geq\ \ 0\ .
\]

This proves the stability of $\mS$.

\end{proof}

\vskip 0.6in


\section{\textbf{Instability of strict graphical strips}}\label{S:proof}

\vskip 0.2in

After these preparations we turn to the core of the proof of Theorem
\ref{T:perminI}. To put the subsequent discussion on a solid ground,
we begin with proving Theorem \ref{T:counterexI}.

\medskip

\begin{proof}[\textbf{Proof of Theorem \ref{T:counterexI}}]
We provide the proof only for the class of surfaces in \eqref{ceI},
leaving it to the interested reader to develop the completely
analogous details for the second class. It is obvious from the
definition \eqref{ceI} that $\mathcal S$ is a $C^2$ graph over the
open subset $\R\times
 I$ of the $(y,t)$-plane. We next observe that $\mathcal S$ has
empty characteristic locus. We can use the global defining function
\begin{equation}\label{df}
\phi(x,y,t)\ =\ x - y G(t)\ , \end{equation} and assume that $\mS$
is oriented in such a way that $\bN = \nabla \phi = (X_1\phi)X_1 +
(X_2\phi) X_2 + (T\phi) T$. Recalling \eqref{pq}, we find
\begin{equation}\label{pqt}
p = X_1\phi = 1 + \frac{y^2}{2} G'(t)\ ,\quad \quad q = X_2 \phi = -
G(t) - \frac{xy}{2} G'(t)\ ,\quad\quad \omega\ = T\phi = - y G'(t)\
.
\end{equation}

Since $p\geq 1>0$, we see from \eqref{pqt} that $\Sigma(\mS) =
\varnothing$. In order to prove the $H$-minimality of $\mathcal S$,
we use \eqref{divY}, which gives
\[
\mathcal H\ =\ X_1 \pb\ +\ X_2 \qb\ .
\]

>From now on, to simplify the notation, we will omit the variable $t$
in all expressions involving $G(t), G'(t), G''(t)$. The second
equation in \eqref{pqt} becomes on $\mS$
\begin{equation}\label{qons}
q\ =\ -\ G\ \left(1\ +\ \frac{y^2}{2} G'\right)\ .
\end{equation}

We thus find on $\mS$
\begin{equation}\label{Wons}
W^2\ =\ p^2 + q^2\ =\ \big(1 + G^2\big) \left(1\ +\ \frac{y^2}{2}
G'\right)^2\ .
\end{equation}

Since they will be useful in the proof of Lemma \ref{L:coefficient},
in what follows we compute several quantities, even if they are not
strictly necessary for the calculation of $\mathcal H$. From
\eqref{pqt} we find
\begin{equation}\label{Xp}
X_1p\ =\ - \frac{y^3}{4} G''\ ,\quad\quad X_2 p\ =\ y G' + \frac{x
y^2}{4} G''\ ,
\end{equation}

\begin{equation}\label{Xq}
X_1q\ =\  \frac{x y^2}{4} G''\ ,\quad\quad X_2 q\ =\ - x G' -
\frac{x^2 y}{4} G''\ .
\end{equation}

>From \eqref{pqt}, \eqref{qons}, \eqref{Xp}, \eqref{Xq} we find on
$\mS$
\begin{equation}\label{Xq2}
X_2 q\ =\ -\ y G \left(G'+ \frac{y^2}{4} G G''\right)\ ,
\end{equation}
and
\begin{equation}\label{X1W}
X_1 W\ =\ \frac{p X_1p + q X_1q}{W}\ =\ - \frac{y^3}{4} G'' \left(1
+ G^2\right)^{\frac{1}{2}}\ , \end{equation}
\begin{equation}\label{X2W}
X_2 W\ =\ \frac{p X_2p + q X_2q}{W}\ =\ y \left(1 +
G^2\right)^{\frac{1}{2}} \left\{G' + \frac{y^2}{4} G G''\right\}\ .
\end{equation}

We now compute
\begin{equation}\label{Hstrips}
\mathcal H\ =\ \frac{X_1p}{W} - p\frac{X_1W}{W^2} + \frac{X_2q}{W} -
q\frac{X_2W}{W^2}\ .
\end{equation}

Using \eqref{pqt}, \eqref{qons}, \eqref{Xp}, \eqref{Xq2},
\eqref{X1W} and \eqref{X2W}, we find
\[
\frac{X_1p}{W} - p\frac{X_1W}{W^2} \ =\ 0\ ,\quad\quad
\frac{X_2q}{W} - q\frac{X_2W}{W^2}\ =\ 0\ .
\]

Inserting the latter two equations in \eqref{Hstrips} we conclude
that $\mS$ is $H$-minimal. To complete the proof of the theorem we
are left with showing that, when $I = \R$, the surface $\mS$ in
\eqref{ceI} is a global intrinsic $X_1$-graph according to
\cite{FSS4}. To prove this, we want to show that there exist
curvilinear coordinates $(u,v)\in \R^2$, and $\phi\in
C^2(\R^2_{u,v})$, such that $\mS$ can be globally parameterized by
\begin{equation}\label{one}
\theta(u,v)\ =\ \left(\phi(u,v),u,v - \frac{u}{2}\phi(u,v)\right)\ .
\end{equation}

We thus see that we must have
\[
\phi(u,v)\ =\ x\ ,\quad u\ = y\ ,\quad\quad v -
\frac{u}{2}\phi(u,v)\ =\ t\ .
\]

These equations give
\[
u = y\ ,\ v = t + \frac{u}{2}\phi(u,v)\ =\ t + \frac{y}{2}x\ =\ t +
\frac{y^2}{2}G(t)\ .
\]

We want to show next that the map $\Phi : \R^2_{y,t} \to \R^2_{u,v}$
given by
\[
\Phi(y,t)\ =\ \left(y, t + \frac{y^2}{2}G(t)\right)\ ,
\]
defines a global diffeomorphism onto. Now, its Jacobian is given by
\begin{equation}\label{Phi}
det\ \begin{pmatrix} 1 & 0
\\
y G(t) & 1 + \frac{y^2}{2}G'(t)
\end{pmatrix}\ =\ 1+ \frac{y^2}{2}G'(t)\ \not= \ 0\ , \ \text{for every}\ (y,t)\in \R^2\ .
\end{equation}

Furthermore, $\Phi$ is globally one-to-one. Assume in fact that
$\Phi(y_1,t_1) = \Phi(y_2,t_2)$, then we have
\begin{equation}\label{two}
y_1 = y_2\ ,\ \quad t_1 + \frac{y_1^2}{2}G(t_1) = t_2 +
\frac{y_2^2}{2}G(t_2)\ .
\end{equation}

Now, let $\alpha = y_1 = y_2$, then either $\alpha= 0$, in which
case \eqref{two} gives $t_1 = t_2$, or $\alpha \neq 0$. In this
second case, we look at the function $f(t) = t + (\alpha^2/2) G(t)$,
and we see that $f$ is strictly increasing over $\R$. Therefore,
\eqref{two} forces again the conclusion $t_1 = t_2$. It is also easy
to see that $\Phi$ is onto. Thanks to the assumption $G'>0$ one has
in fact for every $y\in \R$
\[
\underset{t\to\pm \infty}{\lim}\ t + \frac{y^2}{2} G(t)\ =\ \pm\
\infty\ .
\]

Therefore, given $(u,v)\in \R^2$, if we choose $y=u$, then we can
always find $t\in R$ such that $t + \frac{y^2}{2} G(t)= v$. In
conclusion, $\Phi$ is globally invertible on $\R^2$. Denote by
$\Psi(u,v) = (\Psi_1(u,v),\Psi_2(u,v))$ the inverse of $\Phi$.
Clearly, $y(u,v) = u$. But then the function
\begin{equation}\label{phice1}
\phi(u,v)\ =\ u G(\Psi_2(u,v))\ ,
\end{equation}
defines $\mS$ as a global intrinsic $X_1$-graph. We note that, using
\eqref{Phi} and the inverse function theorem, we obtain for the
Jacobian matrix of $\Psi$
\begin{equation}\label{Psi}
J_\Psi(u,v)\ =\ \begin{pmatrix} 1 & 0
\\
- \frac{u G}{1 + \frac{u^2}{2} G'} & \frac{1}{1 + \frac{u^2}{2} G'}
\end{pmatrix}\ ,
\end{equation}
where for brevity we have written $G$ instead of $G(\Psi_2(u,v))$,
and similarly for $G'$.

In closing, for the benefit of the reader, we provide a second
derivation of the $H$-minimality of the surface \eqref{ceI} based on
the fact that it is locally an intrinsic $X_1$-graph as in
\eqref{one}, with $\phi(u,v)$ given by \eqref{phice1}. We stress
that the following computations only use the fact that $\mS$ be
locally defined as in \eqref{one} in the neighborhood of any fixed
point, for some $C^2$ function $\phi(u,v)$. We will use the
following formula, found in \cite{GS}, \cite{BSV}, for the $H$-mean
curvature of an intrinsic graph in $\HH$
\begin{equation}\label{H}
\mathcal H\ =\ -\ \mathcal B_\phi\left(\frac{\mathcal
B_\phi(\phi)}{\sqrt{1 + \mathcal B_\phi(\phi)^2}}\right)\ ,
\end{equation}
where for a function $f\in C^1(\R)$
\begin{equation}\label{B}
B_\phi(f)\ =\ f_u + \phi f_v
\end{equation}
denotes the linear transport operator. One can easily verify that
\begin{equation}\label{H2}
\mathcal H\ =\ -\ \frac{\mathcal B_\phi(\mathcal
B_\phi(\phi))}{\left(1 + \mathcal
B_\phi(\phi)^2\right)^{\frac{3}{2}}}\ ,
\end{equation}
and therefore the condition that $\mS$ be $H$-minimal becomes
\begin{equation}\label{H3}
\mathcal B_\phi(\mathcal B_\phi(\phi))\ =\ 0\ , \quad\quad\ \phi\in
C^2(\R^2)\ ,
\end{equation}
where now
\begin{equation}\label{burger}
\mathcal B_\phi(\phi)\ =\ \phi_u + \phi \phi_v\ ,
\end{equation}
denotes the nonlinear inviscid Burger operator. Using
\eqref{phice1}, \eqref{Psi}, we compute
\[
\phi_u\ =\ G + u G' \frac{\p \Psi_2}{\p u}\ =\  G + u G' \left(-
\frac{u G}{1 + \frac{u^2}{2} G'}\right)\ =\ \frac{G \left(1 -
\frac{u^2}{2} G'\right)}{1 + \frac{u^2}{2} G'}\ ,
\]
\[
\phi_v\ =\ u G' \frac{\p \Psi_2}{\p v}\ =\ =\ \frac{u G'}{1 +
\frac{u^2}{2} G'}\ .
\]

>From the last two formulas we find
\[
\mathcal B_\phi(\phi)\ =\ G\ .
\]

This gives
\[
\mathcal B_\phi(\mathcal B_\phi(\phi)) \ = \  \mathcal B_\phi(G)\ =\
G_u + \phi G_v\ =\ G' \frac{\p \Psi_2}{\p u} + u G G'\frac{\p
\Psi_2}{\p v}\ = \ - \frac{u G G'}{1 + \frac{u^2}{2} G'} +  \frac{u
G G'}{1 + \frac{u^2}{2} G'}\ =\ 0\ ,
\]
which, according to \eqref{H3}, proves the $H$-minimality of $\mS$.

\end{proof}

\medskip

We now turn to the proof of Theorem \ref{T:perminI}. Since we will
want to compute the second variation of a graphical strip such as
\eqref{ceI} with respect to deformations along the horizontal normal
$\nuX$, we will need to use formula \eqref{2varnu} in Theorem
\ref{T:nudef}. As a first step, we will compute the quantities which
appear as the coefficient of $h^2$ in the integral in the right-hand
side of \eqref{2varnu}. This is the content of the next lemma.

\medskip

\begin{lemma}\label{L:coefficient}
Let $\mS$ be the $H$-minimal surface given by \eqref{ceI}, then one
has
\begin{align}\label{2def}
& 2 (\pb T\qb - \qb T\pb) + 2 \ob(\qb Y\pb - \pb Y\qb) + \ob^2 \ =\
-\ \frac{2 G'(t)}{W^2}\ .
\end{align}
\end{lemma}

\begin{proof}[\textbf{Proof}]
As in the proof of Theorem \ref{T:counterexI}, we use the global
defining function \eqref{df}, and we obtain \eqref{pqt}. Using
\eqref{Y} and \eqref{Xp}, \eqref{Xq}, we obtain on $\mS$
\begin{equation}\label{Yp}
Yp\ =\ \frac{1}{W}\left\{p X_1p + q X_2p\right\}\ =\ -\ \frac{y}{(1
+ G^2)^{\frac{1}{2}}}\left\{\frac{y^2}{4} G''(1 + G^2) +  G
G'\right\}\ ,
\end{equation}

\begin{equation}\label{Yq}
Yq\ =\ \frac{1}{W}\left\{p X_1q + q X_2q\right\}\ =\ \frac{y G}{(1 +
G^2)^{\frac{1}{2}}}\left\{\frac{y^2}{4} G''(1 + G^2) +  G
G'\right\}\ .
\end{equation}

Combining \eqref{Yp}, \eqref{Yq} we find
\begin{equation}\label{YW}
YW\ =\ \frac{1}{W} \left\{p Yp + q Yq\right\}\ =\ -\ y
\left\{\frac{y^2}{4} G''(1 + G^2) +  G G'\right\}\ .
\end{equation}

Combining \eqref{Yp}, \eqref{Yq} and \eqref{YW}, two small miracles
happen, namely
\begin{equation}\label{Ypb}
Y\pb\ =\ \frac{W Yp - p YW}{W^2}\ =\ 0\ ,\ \quad\quad Y\qb\ =\
\frac{W Yq - q YW}{W^2}\ =\ 0\ .
\end{equation}
We now turn to the computation of the derivatives along the
characteristic direction $T$. Differentiating \eqref{pqt} we obtain
on $\mS$
\begin{equation}\label{Tpq}
Tp\ =\ \frac{y^2}{2} G''\ ,\quad\quad\quad Tq\ =\ -\ \left(G' +
\frac{y^2}{2} G G''\right)\ .
\end{equation}

>From \eqref{Tpq} we find
\begin{equation}\label{TW}
TW\ =\ \frac{1}{W} \left\{p Tp + q T q\right\}\ =\ \frac{1}{(1 +
G^2)^{\frac{1}{2}}} \left\{G G' + \frac{y^2}{2} G'' (1 +
G^2)\right\}\ .
\end{equation}

Using \eqref{Tpq} and \eqref{TW} we obtain
\begin{equation}\label{Tpbqb}
T\pb\ =\ -\ \frac{G G'}{W (1 + G^2)}\ , \quad\quad\quad T\qb\ =\ -\
\frac{G'}{W (1 + G^2)}\ .
\end{equation}

>From \eqref{pqt} and \eqref{Tpbqb} we conclude that
\begin{equation}\label{TpTq}
\pb\ T\qb\ -\ \qb\ T\pb \ =\ -\ \frac{G'}{W (1 +
G^2)^{\frac{1}{2}}}\ .
\end{equation}

Finally, combining \eqref{TpTq} with \eqref{Tpq} and \eqref{pqt}, we
obtain
\begin{align*}
& 2 (\pb T\qb - \qb T\pb) + 2 \ob(\qb Y\pb - \pb Y\qb) + \ob^2 \ =\
-\ \frac{2 G'(t)}{W^2}\ ,
\end{align*}
which proves \eqref{2def}.

\end{proof}

\medskip










\medskip

>From Theorem \ref{T:nudef} and Lemma \ref{L:coefficient}, we obtain
the following corollary.

\medskip

\begin{cor}\label{C:X1graph}
Let $\mS$ be an $H$-minimal surface given as in \eqref{ceI}. For any
$h\in C^2_0(\mS)$, one has
\begin{align}\label{2varnugraph}
\mathcal V^H_{II}(\mS; h \nuX)\ & =\ \int_{\mathcal S} |\dels h|^2
d\sigma_H\ -\ 2\  \int_\mS \frac{h^2\ G'(t)}{W^2}\ d\sigma_H\ .
\end{align}
For any $a\in C^2_0(\mS)$, one has
\begin{align}\label{2varX1graph}
\mathcal V^H_{II}(\mS; a X_1)\ & =\ \int_{\mathcal S} \frac{\left(1
+ \frac{y^2}{2} G'(t)\right)^2}{W^2} |\dels a|^2 d\sigma_H\ -\ 2\
\int_\mS \frac{a^2}{W^2 (1 + G(t)^2)} d\sigma_H\ .
\end{align}
\end{cor}

\medskip

To proceed further we next project  onto the $(y,t)$ plane the
formulas \eqref{2varnugraph}, \eqref{2varX1graph} by means of the
$C^2$ parametrization $\theta : \R \times I \to \mathbb R^3$ of the
surface $\mS$ given by $\theta(y,t) = (y G(t),y,t)$.

\medskip

\begin{lemma}\label{L:hardyequiv}
Let $\mS$ be the $H$-minimal surface given by \eqref{ceI}.  For any
$h\in C^2_0(\mS)$, one has
\begin{align}\label{he2}
\mathcal V^H_{II}(\mS; h \nuX)\ & =\ \int_{\Ri} \frac{\left(1 +
\frac{y^2}{2}
G'(t)\right)\ u_y^2}{(1 + G(t)^2)^{1/2}}\ dy dt\\
& -\ 2  \int_{\Ri} \frac{u^2\ G'(t)}{\left(1 + \frac{y^2}{2}
G'(t)\right) (1 + G(t)^2)^{1/2}}\ dy dt\ ,\notag
\end{align}
where we have set $u = h \circ \theta \in C^2_0(\Ri)$. For any $a\in
C^2_0(\mS)$, one has
\begin{align}\label{he}
\mathcal V^H_{II}(\mS; a X_1)\ & =\ \int_{\Ri} \frac{\left(1 +
\frac{y^2}{2}
G'(t)\right)\ u_y^2}{(1 + G(t)^2)^{3/2}}\ dy dt\\
& -\ 2 \int_{\Ri} \frac{u^2\ G'(t)}{\left(1 + \frac{y^2}{2}
G'(t)\right) (1 + G(t)^2)^{3/2}}\ dy dt\ ,\notag
\end{align}
where this time  we have let $u = a \circ \theta \in C^2_0(\Ri)$.
\end{lemma}

\begin{proof}[\textbf{Proof}]
In order to prove \eqref{he2} we make some reductions. Keeping in
mind \eqref{pm}, from \eqref{Wons} we obtain
\begin{equation}\label{hsquare}
\int_\mS \frac{h^2\ G'(t)}{W^2} d\sigma_H\ =\ \int_{\Ri} \frac{u^2\
G'(t)}{\left(1 + \frac{y^2}{2} G'(t)\right) (1 + G(t)^2)^{1/2}}\ dy
dt\ .
\end{equation}

In order to express the first integral in the right-hand side of
\eqref{2varnugraph} as an integral on $\Ri$, we compute $|\dels
h|^2$. We have from \eqref{Y}, \eqref{pqt} and \eqref{qons}
\begin{align}\label{delh}
|\dels h|^2\ & =\ (Zh)^2\ =\ (\qb X_1 h - \pb X_2 h)^2
\\
& =\ \frac{(G(t) X_1 h + X_2 h)^2}{1 + G(t)^2}\ . \notag
\end{align}

Now, the chain rule gives $u_y = G(t)h_x + h_y$, and therefore we
see that we have on $\mS$
\[
G(t) X_1 h + X_2 h\ =\ G(t)\ h_x + h_y\ =\ u_y\ .
\]

>From \eqref{delh} we thus conclude that
\begin{equation}\label{deltasquare}
\int_{\mathcal S} \left|\dels h\right|^2 d\sigma_H\ =\ \int_{\Ri}
\frac{\left(1\ +\ \frac{y^2}{2} G'(t)\right)\ u_y^2}{(1 +
G(t)^2)^{1/2}}\ dy dt\ .
\end{equation}

Combining \eqref{hsquare} and \eqref{deltasquare} we obtain
\eqref{he2}. The proof of \eqref{he} proceeds analogously, and we
omit the details.

\end{proof}

\medskip

The next lemma is the keystone to the proof of Theorem
\ref{T:perminI}. In order to state it, given an interval $\Id =
(-4\delta,4\delta)$, with $\delta>0$, we fix a function $\chi \in
C^\infty_0(\R)$, such that $0\leq \chi(s) \leq 1$, $\chi \equiv 1$
on $|s|\leq \delta$, $\chi \equiv 0$ for $|s|\geq 2\delta$,
$|\chi'|\leq C= C(\delta)$, and $\int_\R \chi(s) ds = 1$. For every
$k\in \mathbb N$ we define $\chi_k(s) = \chi(s/k)$, so that
$\chi_k(s) \equiv 1$ for $|s|\leq \delta k$, $\chi_k(s) \equiv 0$
for $|s| \geq 2\delta k$, and $|\chi_k'(s)| \leq C/k$, with $C$
independent of $k$. We also let $\tilde \chi_k(s) = k \chi(ks)$, and
notice that the $\int_\R \tilde \chi_k(s) ds = 1$ for every $k\in
\mathbb N$, and that $supp(\tilde \chi_k)\subset
[-2\delta/k,2\delta/k]$.

\medskip

\begin{lemma}\label{L:instab}
Let $G\in C^2(\Id)$ be such that $G'>0$ on $\Id$. Define for $k\in
\mathbb N$
\[
f_k(y,t)\ = \ \frac{\chi_{k}(y)\chi(t)}{\sqrt{1 +
\frac{y^2}{2}\,G_k'(t)}}\
\]
where we have let $G_k(t) = G\star \tilde{\chi}_k(t)$. We have
$f_k\in C^\infty_0(\Rd)$, and there exists $k_0 \in \mathbb N$ such
that for all $k > k_0$
\begin{equation}\label{reverse_ineq}
\int_{\Rd} \frac{1 + \frac{y^2}{2}\,G'(t)}{\sqrt{1 + G(t)^2}}\
\left(\frac{\partial f_k}{\partial y}(y,t)\right)^2\,dy dt \ < \ 2
\int_{\Rd} \frac{G'(t)}{(1 + \frac{y^2}{2}\,G'(t))\sqrt{1 +
G(t)^2}}\ f_k(y,t)^2\,dy dt\ .
\end{equation}
\end{lemma}

\begin{proof}[\textbf{Proof}]
We begin by observing that since $G'\in C(\Id)$, such function is
uniformly continuous  on $[-2\delta,2\delta]$. We recall from basic
properties of approximations to the identity that $G_k'= G'\star
\tilde{\chi}_k \to G'$ uniformly on $[-2\delta,2\delta]$. As a
consequence of this, and of the fact that there exists $\epsilon >
0$ such that $G'\geq \epsilon$ on $[-2\delta,2\delta]$, we can find
$k_0\in \mathbb N$ such that for all $k\geq k_0$, and for all $t\in
[-2\delta,2\delta]$, one has
\begin{equation}\label{dominate}
\begin{cases}
\frac{1}{2}\,G'(t)\ \leq\ G_k'(t)\ \leq\ 2\,G'(t)\ ,
\\
\frac{1}{2}\,(1 + \frac{y^2}{2}\,G'(t)) \ \leq \ 1 +
\frac{y^2}{2}\,G_k'(t) \ \leq \ 2\,(1 + \frac{y^2}{2}\,G'(t))\ .
\end{cases}
\end{equation}

We begin with the right-hand side of \eqref{reverse_ineq}.

\begin{align}\label{rhs}
(RHS)\ =\ & 4\ \int_{\R} f_k(y,t)^2\left(\frac{G'(t)}{(2 +
y^2\,G'(t))\sqrt{1 + G(t)^2}}\right)\,dy dt \\
\notag & \qquad \ =\ 8\ \int_{-2\delta}^{2\delta}
\chi(t)^2\,\frac{G'(t)}{\sqrt{1 + G(t)^2}} \left(\int_\R
\frac{\chi_{k}(y)^2}{(2 + y^2\,G'(t))\,(2 +
y^2\,G_k'(t))}\,dy\right)\,dt \\
\notag & \qquad \ \longrightarrow\ 8\ \int_{-2\delta}^{2\delta}
\chi(t)^2\,\frac{G'(t)}{\sqrt{1 + G(t)^2}}\left(\int_\R \frac{1}{(2
+
y^2\,G'(t))^2}\,dy\right)\,dt\\
\notag & \qquad \ = \ 8 \int_{-2\delta}^{2\delta}
\chi(t)^2\,\frac{G'(t)}{\sqrt{1 +
G(t)^2}}\left(\frac{\sqrt{2}\,\pi}{8\,\sqrt{G'(t)}}\right)\,dt \\
\notag &\qquad \ =\ \sqrt{2}\,\pi\int_{-2\delta}^{2\delta}
\chi(t)^2\,\sqrt{\frac{G'(t)}{1 + G(t)^2}}\,dt \qquad\text{as
$k\to\infty$\ ,}
\end{align}
where we have used \eqref{dominate} and Lebesgue dominated
convergence theorem.

On the other hand, we obtain for the integral in the left-hand side
of \eqref{reverse_ineq}
\begin{align}\label{lhs}
& (LHS)\  =\ \int_{\Rd} \frac{W}{1 +
G(t)^2}\left(\frac{\partial}{\partial y}
f_k(y,t)\right)^2\,dy dt \\
\notag &\qquad \ =\ \int_{\Rd} \chi(t)^2 \frac{2 +
y^2\,G'(t)}{\sqrt{1 + G(t)^2}}\, \left\{\frac{\chi_k'(y)^2}{2 +
y^2\,G_k'(t)} - \frac{2\,y\,\chi_k(y)\chi_k'(y)\,G_k'(t)}{(2 +
y^2\,G_k'(t))^2} + \frac{y^2\,\chi_k(y)^2\,G_k'(t)^2}{(2 +
y^2\,G_k'(t))^3}\right\}\,dydt
\\
\notag &\qquad \ =\ \int_{\Rd}\frac{\chi(t)^2}{\sqrt{1 + G(t)^2}}\
\chi_k'(y)^2\ \frac{2
+ y^2\,G'(t)}{2+y^2\,G_k'(t)}\,dydt \\
\notag &\qquad\qquad \ -\
\int_{\Rd}\chi(t)^2\,\frac{G_k'(t)}{\sqrt{1 + G(t)^2}} \frac{y\,(2 +
y^2\,G'(t))}{(2 + y^2\,G_k'(t))^2}
\bigl(\chi_k(y)^2\bigr)^{'} \,dydt \\
\notag &\qquad\qquad \ +\
\int_{\Rd}\chi(t)^2\,\frac{G_k'(t)^2}{\sqrt{1 + G(t)^2}}\
\chi_k(y)^2\ \frac{y^2\,(2 + y^2\,G'(t))}{(2 +
y^2\,G_k'(t))^3}\,dydt\ .
\end{align}

Using the fact that
\[
\frac{\partial}{\partial y}\left(\frac{y}{2 + y^2\,G_k'(t)}\right) \
=\ \frac{2 - y^2\,G_k'(t)}{(2 + y^2\,G_k'(t))^2}\ , \qquad\quad
\frac{\partial}{\partial y}\left(\frac{2 + y^2\,G'(t)}{2 +
y^2\,G_k'(t)}\right) \ \ = \ \frac{4\,y\,(G'(t) - G_k'(t))}{(2 +
y^2\,G_k'(t))^2}\ ,
\]
we integrate by parts the integral containing the term
$\bigl(\chi_k(y)^2\bigr)^{'}$, obtaining
\begin{align}\label{secondTerm}
& \ -\ \int_{\Rd}\chi(t)^2\,\frac{G_k'(t)}{\sqrt{1 + G(t)^2}}
\frac{y\,(2 + y^2\,G'(t))}{(2 + y^2\,G_k'(t))^2}
\bigl(\chi_k(y)^2\bigr)^{'} \,dydt \\
\notag &\qquad \ =\ 2\,\int_{\Id} \chi(t)^2\frac{G_k'(t)}{\sqrt{1 +
G(t)^2}}\left(\int_\R \frac{\chi_k(y)^2}{(2 +
y^2\,G_k(t)^2)^2}\,\frac{2 + y^2\,G'(t)}{2 +
y^2\,G_k'(t)}\,dy\right)\,dt \\
\notag &\qquad\qquad \ -\
\int_{\Rd}\chi(t)^2\,\frac{G_k'(t)^2}{\sqrt{1 + G(t)^2}}\
\chi_k(y)^2\
\frac{y^2\,(2 + y^2\,G'(t))}{(2 + y^2\,G_k'(t))^3}\,dydt \\
\notag &\qquad\qquad \ +\ \int_{\Id} \chi(t)^2\,\frac{G_k'(t)(G'(t)
- G_k'(t))}{\sqrt{1 + G(t)^2}}\, \left(\int_\R
\chi_k(y)^2\,\frac{4\,y^2}{(2 + y^2\,G_k'(t))^2}\,dy\right)\,dt\ .
\end{align}

Using \eqref{secondTerm} in \eqref{lhs}, and after canceling one
term, we obtain for the left-hand side of \eqref{reverse_ineq}
\begin{align}\label{lhs2}
& (LHS)\ =\ \int_{\Rd} \frac{W}{1 +
G(t)^2}\left(\frac{\partial}{\partial y}
f_k(y,t)\right)^2\,dy dt \\
\notag &\qquad \ =\ \int_{\Rd}\frac{\chi(t)^2}{\sqrt{1 + G(t)^2}}\
\chi_k'(y)^2\ \frac{2
+ y^2\,G'(t)}{2+y^2\,G_k'(t)}\,dydt \\
\notag & \qquad\qquad \ +\ 2\,\int_{-2\delta}^{2\delta}
\chi(t)^2\frac{G_k'(t)}{\sqrt{1 + G(t)^2}}\,\left(\int_\R
\frac{\chi_k(y)^2}{(2 + y^2\,G_k'(t))^2}\,\frac{2 + y^2\,G'(t)}{2 +
y^2\,G_k'(t)}\,dy\right)\,dt\\
\notag & \qquad\qquad \ +\ \int_{\Id} \chi(t)^2\,\frac{G_k'(t)(G'(t)
- G_k'(t))}{\sqrt{1 + G(t)^2}}\, \left(\int_\R
\chi_k(y)^2\,\frac{4\,y^2}{(2 + y^2\,G_k'(t))^2}\,dy\right)\,dt
\notag\\
& =\ I_k\ +\ II_k\ +\ III_k\ . \notag
\end{align}

We now analyze each of the three integrals in \eqref{lhs2}.  Using
the fact that $G'(t)>0$ on $\Id$, we also have $G_k'(t) = G'\star
\tilde{\chi}_k(t)
> 0$ on $\Id$. When $k\to \infty$ the first integral satisfies
\begin{align}\label{firsTerm}
I_k\ =\ & \int_{\Ri}\frac{\chi(t)^2}{\sqrt{1 + G(t)^2}}\
\chi_k'(y)^2\
\frac{2 + y^2\,G'(t)}{2+y^2\,G_k'(t)}\,dydt \\
\notag &\qquad \leq\ \frac{C}{k^2}\,\int_{-2\delta}^{2\delta}
\frac{1}{\sqrt{1 + G(t)^2}}\left(\int_{-2k\delta}^{2k\delta}\frac{2
+ y^2\,G'(t)}{2 +
y^2\,G_k'(t)}\,dy\right)\,dt \\
\notag \text{(by \eqref{dominate})} &\qquad \leq\ \frac{C}{k^2}
\int_{-2\delta}^{2\delta} \frac{1}{\sqrt{1 +
G'(t)^2}}\,\left(\int_{-2k\delta}^{2k\delta} 2\,dy\right)\,dt \ =\
\frac{8C\delta}{k}\,\int_{-2\delta}^{2\delta} \frac{1}{\sqrt{1 +
G'(t)^2}}\,dt \longrightarrow\ 0\ .
\end{align}

Similarly, letting $k\to\infty$ we have for the third integral
\begin{align}\label{thirdTerm}
III_k\ =\ & \int_{\Id} \chi(t)^2\,\frac{G_k'(t)(G'(t) -
G_k'(t))}{\sqrt{1 + G(t)^2}}\,
\left(\int_\R \chi_k(y)^2\,\frac{4\,y^2}{(2 + y^2\,G_k'(t))^2}\,dy\right)\,dt\\
\notag \text{(by \eqref{dominate})} & \quad\quad \leq\
16\,\underset{t\in[-2\delta,2\delta]}{sup} \bigl|G_k'(t) -
G'(t)\bigr|\,\int_{-2\delta}^{2\delta} \frac{G'(t)}{\sqrt{1 +
G(t)^2}}\left(\int_\R \frac{y^2}{(2 + y^2\,G'(t))^2}\,dy\right)\,dt
\longrightarrow\ 0\ ,
\end{align}
from the uniform convergence of $G_k'$ to $G'$ on
$[-2\delta,2\delta]$. Finally, we have
\begin{align}\label{goodTerm}
II_k \ =\ & 2\,\int_{-2\delta}^{2\delta}
\chi(t)^2\frac{G_k'(t)}{\sqrt{1 + G(t)^2}}\,\left(\int_\R
\frac{\chi_k(y)^2}{(2 + y^2\,G_k'(t))^2}\,\frac{2 + y^2\,G'(t)}{2 +
y^2\,G_k'(t)}\,dy\right)\,dt \\
\notag &\qquad \longrightarrow\ 2\ \int_{-2\delta}^{2\delta}
\chi(t)^2\, \frac{G'(t)}{\sqrt{1 +
G(t)^2}}\left(\int_\R \frac{1}{(2 + y^2\,G'(t)^2)^2}\,dy\right)\,dt \\
\notag &\qquad\qquad \ =\ 2 \int_{-2\delta}^{2\delta}
\chi(t)^2\,\frac{G'(t)}{\sqrt{1 + G(t)^2}}
\left(\frac{\sqrt{2}\,\pi}{8\,\sqrt{G'(t)}}\right)\,dt \ =\
\frac{\sqrt{2}\,\pi}{4}\ \int_{-2\delta}^{2\delta} \chi(t)^2\,
\sqrt{\frac{G'(t)}{1 + G(t)^2}}\,dt\ .
\end{align}
In the above, we have used \eqref{dominate} to deduce that for large
enough $k$
\[
\frac{G_k'(t)}{\sqrt{1 + G(t)^2}}\left(\frac{\chi_k(y)^2}{(2 +
y^2\,G_k'(t))^2}\,\frac{2 + y^2\,G'(t)}{2 + y^2\,G_k'(t)}\right) \
\leq\ 16\,\frac{G'(t)}{\sqrt{1 + G(t)^2}}\, \frac{1}{(2 +
y^2\,G'(t))^2}\in L^1(\R\times [-2\delta,2\delta])
\]
and therefore, Lebesgue dominated convergence theorem applies. To
summarize, we have as $k\to\infty$
\[
(LHS)\ =\ \int_{\Rd} \frac{W}{1 +
G(t)^2}\left(\frac{\partial}{\partial y} f_k(y,t)\right)^2\,dy dt
\longrightarrow \frac{\sqrt{2}\,\pi}{4}\ \int_{-2\delta}^{2\delta}
\chi(t)^2\,\sqrt{\frac{G'(t)}{1 + G(t)^2}}\,dt\ .
\]
Combining this with\eqref{rhs}, we reach the sought for conclusion.

\end{proof}

\medskip

We are now ready to prove Theorem \ref{T:perminI}.

\begin{proof}[\textbf{Proof of Theorem \ref{T:perminI}}]
Let $\tilde \mS$ be a graphical strip, then there exist $I\subset
\R$, $\tilde G\in C^2(\R)$, with $G\geq 0$, such that, after
possibly a left-translation and a rotation about the $t$-axis,
$\tilde \mS$ can be represented in the form $x= y \tilde G(t)$ for
$(y,t)\in \R\times I$. If we assume further that $\tilde \mS$ is a
strict graphical strip, then we can find an interval $J =
(a,b)\subset I$, such that $\tilde G'>0$ on $J$. Since the
stability, or the instability, are invariant under left-translations
and rotations, it will suffice to prove that $\tilde \mS$ is
unstable. Assume without restriction that $- \infty < a < b <
\infty$, and set $t_0 = (a+b)/2$, $g_0= (0,0,-t_0)$. Consider the
left-translated surface $\mS = g_0\circ \tilde \mS$, see \eqref{Hn},
then $\mS$ is described by
\[
x\ =\ y\  G(t)\ , \quad\quad\quad (y,t)\in \R \times \Id\ ,
\]
with $4 \delta = (b-a)/2$, and $G(t) = \tilde G(t_0 + t)$. Since it
is clear that $G'>0$ on $\Id$, we can apply Lemma \ref{L:instab},
and conclude that there exists $k_0\in \mathbb N$ such that for
$k\geq k_0$ the sequence $f_{k}$ satisfies \eqref{reverse_ineq}.
This being said, we now define $h_k:\HH\to\R$ as follows
\[
h_k(x,y,t)\ =\ \frac{1}{\sqrt{1 +
\frac{y^2}{2}G_k'(t)}}\,\chi_{k}(y)\chi(t)\chi_{k}(x - y\,G(t))\ .
\]

We observe that $h_k(\theta(y,t)) = f_k(y,t) \chi_k(0) = f_k(y,t)$,
and therefore $h_k\in C^2_0(\mS)$. At this point, appealing to
\eqref{he2} in Lemma \ref{L:hardyequiv}, and to Lemma
\ref{L:instab}, we conclude that for every fixed $k\geq k_0$, we
have
\[
\mathcal V^H_{II}(\mS;h_k \nuX)\ <\ 0\ .
\]

This proves that $\mS$ is unstable, and therefore such is also the
surface $\tilde \mS$.

\end{proof}

\medskip

\begin{rmrk}\label{R:global}
In particular, since every global minimizer is also a local one, we
have also shown that $\mS$ cannot be a global minimizer of the
$H$-perimeter.
\end{rmrk}

\vskip 0.6in


\section{\textbf{Instability of $H$-minimal entire graphs and proof of the Bernstein conjecture}}\label{S:rs}

\vskip 0.2in

In this section we prove Theorems \ref{T:reductionI} and
\ref{T:instability}. Our strategy will be to first establish Theorem
\ref{T:reductionI}, and then combine this result with Theorem
\ref{T:perminI} to obtain Theorem \ref{T:instability}. As we have
mentioned in the introduction, our proof of Theorem
\ref{T:reductionI} is based on the main results in \cite{GP}, but we
reiterate that an alternative proof could be obtained combining the
results in the two papers \cite{CHMY} and \cite{CH}.

We begin by recalling the basic notion of seed curve from \cite{GP}.
In what follows, $\Om \subset \R^2$ denotes a given connected, open
set of the $(x,y)$-plane, and $f\in C^k(\Om)$, with $k\geq 2$. We
consider the graph of $f$ over $\Om$
\begin{equation}\label{xygraph}
\mS\ =\ \{(x,y,t)\in \HH\mid (x,y) \in \Omega , t = f(x,y)\}\ .
\end{equation}

We assume that $\mS$ be oriented in such a way that $\bN = \nabla
\phi = (X_1\phi) X_1 + (X_2\phi) X_2 + (T\phi) T$, were $\phi(x,y,t)
= t - f(x,y)$. From \eqref{pq} we obtain
\[
p\ =\ X_1\phi\ =\ - f_x - \frac{y}{2}\ , \quad q = X_2 \phi\ =\ -
f_y + \frac{x}{2}\ ,\quad\quad W = \sqrt{\left(f_x +
\frac{y}{2}\right)^2 + \left(f_y - \frac{x}{2}\right)^2}\ ,
\]
and thus we have from \eqref{normals},
\[
\nuX\ =\ -\ \frac{f_x + \frac{y}{2}}{\sqrt{\left(f_x +
\frac{y}{2}\right)^2 + \left(f_y - \frac{x}{2}\right)^2}}\ X_1\ -\
\frac{f_y - \frac{x}{2}}{\sqrt{\left(f_x + \frac{y}{2}\right)^2 +
\left(f_y - \frac{x}{2}\right)^2}}\ X_2\ ,
\]
away from $\Sigma(\mS)$. We stress that, since $\mS$ is a graph over
the $(x,y)$-plane, $\nuX$ is independent of the variable $t$. Such
crucial property would not be true for a graphical portion over the
coordinate planes $(y,t)$ or $(x,t)$. We can thus identify in a
natural fashion $\nuX$ with a unit $C^{k-1}$ vector field $\til$
onto the $(x,y)$-plane as follows (recall that $k\geq 2$). Given a
point $g = (z,t)\in \mS$, with $z=(x,y)$, we let
\[
\til(z)\ =\  \left(-\ \frac{f_x + \frac{y}{2}}{\sqrt{\left(f_x +
\frac{y}{2}\right)^2 + \left(f_y - \frac{x}{2}\right)^2}}\ ,\ -\
\frac{f_y - \frac{x}{2}}{\sqrt{\left(f_x + \frac{y}{2}\right)^2 +
\left(f_y - \frac{x}{2}\right)^2}}\right)\ . \]

The notation
\[
(\til)^\perp\ =\ \left(-\ \frac{f_y - \frac{x}{2}}{\sqrt{\left(f_x +
\frac{y}{2}\right)^2 + \left(f_y - \frac{x}{2}\right)^2}}\ ,\
\frac{f_x + \frac{y}{2}}{\sqrt{\left(f_x + \frac{y}{2}\right)^2 +
\left(f_y - \frac{x}{2}\right)^2}}\right)\ ,
\]
will indicate the unit vector field in $\Om$ perpendicular to $\til$
(with respect to the Euclidean inner product $u\cdot v$ in $\R^2$).

\medskip

\begin{dfn}\label{D:seeddef}  Let $\mS$ be a $C^2$ graph as in
\eqref{xygraph}, with $\Sigma(\mS) = \varnothing$, and suppose that
$\mS$ be $H$-minimal. Given a point $z\in \Om\subset \R^2$, a
\emph{seed curve} of $\mS$ based at $z$ is defined to be the
integral curve of the vector field $\til$ with initial point $z$.
Denoting such a seed curve by $\gamma_{z}(s)$, we then have
  \begin{equation}\label{D:seedcurve}
\gamma'_z(s)\ =\
  \til(\gamma_z(s))\ ,\quad\quad\quad \gamma_z(0)\ =\ z\ .
  \end{equation}

We will indicate by $\tilL$ the integral curve of $(\til)^\perp$
starting at the point $z$, i.e.,
\begin{equation}\label{seedp} \tilLp\ =\ \tilp(\tilL)\
,\quad\quad\quad \tilLzero\ =\ z\ .
\end{equation}
 If the base point $z$ is understood or irrelevant, we simply denote the seed curve by $\gamma(s)$, and $\tilL$ by $\tilde{\mathcal{L}}(r)$.
\end{dfn}

\medskip

We emphasize that the assumption of $H$-minimality on $\mS$, implies
the crucial property that the vector field $\til$ be divergence free
in $\Om$, with respect to the standard divergence operator in
$\R^2$. As a consequence of this fact, it is proved in \cite{GP}
that the integral curves $\tilde{\mathcal{L}}_z$ of $(\til)^\perp$
are straight-line segments. Furthermore, since $|\til|\equiv 1$ in
$\Om$, every seed curve is parameterized by arc-length, and it is a
$C^1$ embedded curve in $\R^2$ over its interval of definition. The
curves $\{\tilde{\mathcal{L}}_z, \gamma_z\}$ are used in \cite{GP}
to define a local $C^1$ diffeomorphism of the $(x,y)$-plane given by
\begin{equation}\label{newpar}
(s,r)\ \ra\ (x(s,r), y(s,r))\ \overset{def}{=}\ F(s,r)\ =\
\gamma(s)\ +\ r\;\tilp(\gamma(s))\ .
\end{equation}

We note explicitly that $F$ maps the straight line $r=0$ into the
seed curve $\gamma(s)$, i.e., \[ F(s,0)\ =\ \gamma(s)\ .
\]

On the other hand, the straight line $s=0$ is mapped into the
straight line passing through the base point $z$ of the seed curve
and having direction vector $\tilp(z)$, i.e., \[ F(0,r)\ =\ z\ +\ r\
\tilp(z)\ . \]

One recognizes that $F(0,r) = \tilL$. In particular, when $z =
\gamma(s)$ we obtain from this identity
\[
\tilLg\ =\ \gamma(s)\ +\ r\ \tilp(\gamma(s))\ .
\]

Along the seed curve we have $\tilp(\gamma(s)) = \gamma'(s)^\perp =
(\gamma_2'(s) , - \gamma_1'(s))$. One thus has the following
explicit expression for $F(s,r)$
\begin{equation}\label{explicit}
F(s,r)\ =\ \gamma(s) + r \gamma'(s)^\perp\ =\ \bigg(\gamma_1(s) + r
\gamma_2'(s)\ ,\ \gamma_2(s) - r \gamma_1'(s)\bigg)\ .
\end{equation}

Henceforth, given $f\in C^k(\Om)$ as in \eqref{xygraph}, we will use
the notation \[ h(s,r)\ = \ f(F(s,r))\ =\ f(\gamma_1(s) + r
\gamma_2'(s)\ ,\ \gamma_2(s) - r \gamma_1'(s))\ , \] for all values
$(s,r)$ for which the right-hand side is defined, and also let
$h_0(s) = h(s,0)$. Notice that \[ h_0(0)\ =\ h(0,0)\ =\
f(\gamma(0))\ =\ f(z)\ . \]

We next define
\begin{equation}\label{para}
(s,r)\ \longrightarrow\ \mathcal F(s,r)\ =\
(\gamma_1(s)+r\gamma_2'(s),\gamma_2(s)-r\gamma_1'(s), h(s,r))\ ,
\end{equation}
and observe explicitly that $\mathcal F(0,0) = (\gamma(0),h(0,0)) =
(z,f(z))$. We will need the following result, which is Theorem A in
\cite{GP}. This result shows, in particular, that, thanks to the
assumption that $\mS$ be $H$-minimal, the function $h(s,r)$ must
take up a special structure, see \eqref{para2} below.

\medskip

\begin{thrm}\label{T:thmA}
Let $\mS\subset \HH$ be a $C^k$ graph of the type \eqref{xygraph},
with $\Sigma(\mS)=\varnothing$. If $\mS$ is $H$-minimal, then for
every $g = (z,t) \in \mS$, there exist  intervals $I,J$, and an open
neighborhood of $g$ on $\mS$ that can be parameterized by
\eqref{para}, where $h(s,r)$ is given by
\begin{equation}\label{para2}
h(s,r)\ =\ h_0(s)\ -\ \frac{r}{2}\ \gamma(s)\cdot \gamma'(s)\
,\quad\quad(s,r)\in I\times J\ ,
\end{equation}
with
\[\gamma \in C^{k+1}(I),\;  h_0 \in C^k(I)\ .
\]
\end{thrm}

\medskip

\begin{cor}\label{C:seednux}  Let $\mS\subset \HH$ be given as in Theorem \ref{T:thmA}. At every $g= \mathcal F(s,r)\in \mS$, one has
\begin{equation}\label{W}
W\ =\ \left|h_0' - r + \frac{r^2}{2} \kappa + \frac{1}{2} \gamma'
\cdot \gamma^\perp\right|\ ,
\end{equation}
where $W$ is the angle function defined in \eqref{pq}. Moreover, the
horizontal Gauss map is given by
\begin{equation}\label{hgmsr}
\nuX\ =\ sgn \left(h_0' - r + \frac{r^2}{2} \kappa + \frac{1}{2}
\gamma' \cdot \gamma^\perp\right) \bigg\{\gamma'_1(s) \; X_1 +
\gamma_2'(s) \; X_2\bigg\}\ .
\end{equation}
\end{cor}

\begin{proof}[\textbf{Proof}]
Suppose that $\theta : \Om \subset \R^2 \to \HH$, with $\theta(u,v)
= x(u,v) X_1 + y(u,v) X_2 + t(u,v) T$ be a local parameterization os
$\mS$, then a direct calculation shows that the coefficients of $\bN
= \theta_u \wedge \theta_v$, and $\bN^H$ as in \eqref{normals}, with
respect to the orthonormal basis $\{X_1,X_2,T\}$, are given by the
equations
\begin{equation}\label{ppar}
\begin{cases}
p\  =\ y_u t_v - y_v t_u\ -\ \frac{y}{2} (x_u y_v - x_v y_u)\ ,
\\
q\  =\ x_v t_u - x_u t_v\ +\ \frac{x}{2} (x_u y_v - x_v y_u)\ ,
\\
\omega\ =\ x_u y_v - x_v y_u\ . \end{cases}
\end{equation}

Applying the formulas \eqref{ppar} to the parameterization
$\theta(u,v) = \mathcal F(s,r)$ given by \eqref{para},
\eqref{para2}, we find
\begin{equation}\label{ppar2}
\begin{cases}
p\  =\ \gamma_1'\left(h_0' - r +\frac{r^2}{2}\kappa + \frac{1}{2}
\gamma' \cdot \gamma^\perp\right)\  ,
\\
q\  =\ \gamma_2'\left(h_0' - r +\frac{r^2}{2}\kappa + \frac{1}{2}
\gamma' \cdot \gamma^\perp\right)\  ,
\\
\omega\ =\ -\ (1 - r \kappa)\ , \end{cases}
\end{equation}
where, following \cite{GP}, we have denoted by
\[
\kappa(s)\ =\ \gamma''(s)\cdot \gamma'(s)^\perp
\]
the signed curvature of the seed curve $\gamma(s)$. From
\eqref{ppar2} we obtain
\[
W\ =\ \sqrt{p^2 + q^2}\ =\ \left|h_0' - r + \frac{r^2}{2} \kappa +
\frac{1}{2} \gamma' \cdot \gamma^\perp\right|\ .
\]

The assumption that the characteristic locus of $\mS$ be empty
implies that for every $(s,r)$ in the domain of $W(s,r)$, one has
\[
h_0' - r + \frac{r^2}{2} \kappa + \frac{1}{2} \gamma' \cdot
\gamma^\perp\ \neq\  0\ .
\]

>From the expression of $W$, and from \eqref{ppar2}, we thus conclude
that the horizontal Gauss map of $\mS$ is given by \eqref{hgmsr}.

\end{proof}

\medskip

Theorem \ref{T:thmA} admits the following converse, which is Theorem
B in \cite{GP}.

\medskip

\begin{thrm}\label{T:thrmAconv}
Given an interval $I\subset \R$, $k \ge 2$, a properly embedded
plane curve $\gamma \in C^{k}(I)$, parameterized by arc-length, and
a function $h_0\in C^{k}(I)$, let $\mS\subset \HH$ be the surface
parameterized by $\mathcal F:I\times \R$ as in \eqref{para}, with
$h(s,r)$ given by \eqref{para2}. Then, $\mS$ is a $C^{k-1}$
$H$-minimal surface.
\end{thrm}

\medskip

Combining Theorems \ref{T:thmA} and \ref{T:thrmAconv}, we conclude
that to specify a patch of a smooth $H$-minimal surface, one must
specify a single curve in $\HH$ determined by a seed curve $\gamma$,
parameterized by arc-length, and an initial height function $h_0$.
We will also need the following result, which is either Theorem A in
\cite{CHMY}, or Theorem E in \cite{GP}.

\medskip

\begin{thrm}\label{T:xygraphs}
Suppose that $\mS\subset \HH$ be a connected $H$-minimal entire
graph over the $(x,y)$-plane, then:
\begin{enumerate}
\item Either $\mS$ is a plane of the form $ax+by + ct = \gamma$ for some real numbers
  $a,b, c, \gamma$, with $c\not= 0$.
\item Or, there exist $g_0=(x_0,y_0,t_0)\in \HH$, $a,b\in \R$ such that $a^2+b^2=1$, and $h_0\in C^2(\R)$, such that $\mS$
is globally parameterized by
\[\left (x+x_0,y+y_0,t_0-\frac{1}{2}ab(x^2-y^2)-\frac{1}{2}(b^2-a^2)xy+h_0(ax+by)+\frac{1}{2}x_0y-\frac{1}{2}xy_0
\right)\ .
\]
\end{enumerate}
\end{thrm}

\medskip

We emphasize that both types of surfaces arising in this theorem
have non-empty characteristic loci. For instance, in the case (1) we
have $\Sigma(\mS) = \{(-2b/c,2a/c,\gamma/c))\}$. Finally, we recall
a basic result in $H$-minimal surface theory. We mention that the
next result is one half of Theorem C in \cite{GP}, and of Corollary
4.2 in \cite{CHMY}.

\medskip

\begin{thrm}\label{foliation}  Let $S$ be a $C^2$, $H$-minimal
surface without boundary, which is complete, connected and embedded,
then $S$ is a ruled surface, all of whose rules are horizontal
straight lines, which are the integral curves of $(\nuX)^\perp$.
\end{thrm}

\medskip

After these preliminaries, we turn to the proof of Theorem
\ref{T:reductionI}. To prepare for it, we first show that an
$H$-minimal surface satisfying the hypothesis of Theorem
\ref{T:reductionI}, can be reduced to an $H$-minimal graph over the
$(y,t)$-plane, having similar properties. We will need some simple
lemmas which clarify the effect of left-translations on a graph. As
we have noted in the introduction, the left-translations \eqref{Hn}
are affine transformations, thereby they preserve planes and lines,
see \eqref{ltplaneI}. Moreover, the left-translations preserve the
property of a surface of having empty characteristic locus, they
preserve the $H$-mean curvature, and therefore the $H$-minimality,
the $H$-perimeter, and the property of a surface of being stable or
unstable. We note that rotations about the $t$-axis (the group
center), also have the same properties.

\medskip

\begin{lemma}\label{L:translateplane}
Suppose $P \subset \HH$ be a plane with Euclidean normal of the form
$\bN^e_P = (a,0,c)$, with $a^2 + c^2\not=0$, and let $\mS$ be a
graph over $P$. For a fixed $g_0=(x_0,y_0,t_0) \in \HH$, $g_0 \circ
\mS$ is a graph over a plane $\tilde P$ with Euclidean normal vector
given by $\bN^e_{\tilde P} = \left (a,0,c-\frac{a y_0}{2}\right)$.
\end{lemma}

\begin{proof}[\textbf{Proof}]
For $g=(x,y,t) \in P$, consider the straight-line through $g$ and
parallel (with respect to the Euclidean inner product in $\R^3$) to
$\bN^e_P$, $L(g) = \{\ell_g(s)=g+s(a,0,c)\in \HH\mid s\in \R\}$. The
assumption that $\mS$ be a graph over $P$ implies the existence of a
unique $s_0\in \R$ such that $L(g) \cap \mS = \{\ell_g(s_0)\}$.
Consider the left-translated line $g_0\circ L(g)$, which is given by
\[
g_0\circ l_g(s)\ =\ g_0\circ g +s\left (a,0,c -
\frac{ay_0}{2}\right)\ ,\quad\quad\quad s\in \R\ ,
\]
and note that $g_0\circ l_g(s_0) \in g_0 \circ \mS$. We see that
$g_0 \circ \mS$ is a graph over a plane $\tilde P$ with Euclidean
normal $\left(a,0,c-\frac{ay_0}{2}\right)$, unless the Euclidean
normal to the plane $g_0 \circ P$ is perpendicular (with respect to
the standard inner product in $\R^3$) to $\left
(a,0,c-\frac{ay_0}{2}\right)$. But this cannot happen. To verify
this, observe that from \eqref{ltplaneI} the Euclidean normal to
$g_0\circ P$ is given by $\left(a+ \frac{cy_0}{2},- \frac{cx_0}{2},
c \right)$. Since the Euclidean inner product of this vector with
$\left(a,0,c-\frac{y_0a}{2}\right)$ is $a^2+c^2 \neq 0$, we reach
the desired conclusion.

\end{proof}

\medskip

We are now ready to accomplish our first reduction.

\medskip

\begin{lemma}\label{L:ytplane}
Let $\mS\subset \HH$ be an $H$-minimal entire graph, with
$\Sigma(\mS) = \varnothing$, and assume that $\mS$ is not itself a
vertical plane such as \eqref{vp}. After composing with a suitable
rotation about the $t$-axis and with a left-translation, we may
assume that there exist $\psi\in C^2(\R^2)$ for which
\[
\mS\ =\ \{(x,y,t)\in \HH\mid (y,t)\in \R^2 , x = \psi(y,t)\}\ .
\]
\end{lemma}

\begin{proof}[\textbf{Proof}]

Suppose that $\mS$ be a graph over the plane $P$ given by
$ax+by+ct=\gamma$ for $a,b,c,\gamma \in \R$, with $a^2 + b^2 + c^2
\neq 0$. Suppose first that $b = 0$, and consider the two cases,
$a=0$, and $a \neq 0$. In the former case, keeping in mind that
$c\neq0$, we see that $\mS$ is a graph over the plane $t=\gamma/c$.
The left-translation by $g_0 = (0,0,- \gamma/c)$, sends this to the
plane $t = 0$, and the surface $g_0\circ \mS$ becomes a global
$H$-minimal graph, with empty characteristic locus, over the
$(x,y)$-plane. This however contradicts Theorem \ref{T:xygraphs},
since this result forces such a graph to have non-empty
characteristic locus. We must thus have $a\neq 0$. In this case,
left-translating by $g_0=(0,-2c/a,0)$, Lemma \ref{L:translateplane}
shows that $g_0 \circ \mS$ is an entire graph over the
$(y,t)$-plane. Thus, we can find $\psi\in C^2(\R^2)$ which defines
$g_0 \circ \mS$, thus yielding the desired conclusion.

We are thus left with the situation $b\neq 0$. In this case,
however, performing a rotation of an angle $\theta = cotan^{-1}(-
a/b)\in (0,\pi)$ about the $t$-axis, which preserves $H$-minimality,
we obtain a new surface which is an entire graph over a plane having
Euclidean normal in the form $(a,0,c)$, and we can thus argue as in
the first part to reach the sought for conclusion.

\end{proof}

\medskip

As it will be useful later, we next prove that the notions of seed
curve and height function are preserved under left-translation.

\medskip

\begin{lemma}\label{L:translateseed}
Suppose there exist intervals $I, J\subset \R$ so that a portion,
$\mS_0$, of an $H$-minimal surface, $\mS$, is parameterized by a
seed curve, $\gamma$, and height function, $h_0$, as in \eqref{para}
and \eqref{para2}, with $s
  \in I, r \in J$. If $g_0=(x_0,y_0,t_0) \in \HH$,
then the surface $g_0 \circ \mS_0$ is also parameterized as in
\eqref{para} and \eqref{para2}, using a seed curve, $\Hat{\gamma}$,
and height function, $\Hat{h}_0$, given by $\Hat{\gamma}(s)=(x_0 +
\gamma_1(s),y_0 + \gamma_2(s))$, and
$\Hat{h}_0(s)=h_0(s)+\frac{x_0}{2}\gamma_2(s)-\frac{y_0}{2}\gamma_1(s)$,
for $s \in I, r \in J$.
\end{lemma}

\begin{proof}[\textbf{Proof}]
With $F(s,r)$ as in \eqref{explicit}, consider  $g_0 \circ \mS_0$,
which is parameterized by $g_0 \circ \mathcal F(s,r) = g_0\circ
(F(s,r),h(s,r))$. Using \eqref{Hn}, we see that $g_0 \circ \mathcal
F(s,r)$ is given by
\[
\Hat{\mathcal F}(s,r)\ =\ (\Hat{F}(s,r),\Hat{h}(s,r))\ ,
\]
with
\[\Hat{F}(s,r)=\Hat{\gamma}(s)+r
(\Hat{\gamma}'(s))^\perp =\gamma(s)+r(\gamma'(s))^\perp + (x_0,y_0)\
,
\]
and
\[\Hat{h}(s,r) = \Hat{h}_0(s)-\frac{r}{2} \Hat{\gamma}(s) \cdot
\Hat{\gamma}'(s)\] where $\Hat{h}_0(s) =
h_0(s)+\frac{x_0}{2}\gamma_2(s)-\frac{y_0}{2}\gamma_1(s)$. Applying
Corollary \ref{C:seednux} to the parameterization $\Hat{\mathcal
F}(s,r)$, we find $\Hat{W}(s,r) = W(s,r)$, and therefore from
\eqref{hgmsr} we obtain for horizontal Gauss map of $g_0 \circ
\mS_0$ at the point $g_0 \circ (F(s,r),h(s,r))$,
\[ sgn \left(h_0' - r + \frac{r^2}{2} \kappa + \frac{1}{2}
\gamma' \cdot \gamma^\perp\right) \bigg\{\Hat{\gamma}_1'(s)\; X_1 +
\Hat{\gamma}_2'(s) \; X_2\bigg\}\ . \]

Since $\Hat{\gamma}_i'(s)=\gamma_i'(s)$, we have that the components
of the horizontal Gauss map are the same as those of the unit
horizontal Gauss map for $\mS_0$ at the point $(F(s,r),h(s,r))$, see
Corollary \ref{C:seednux}. Upon projection to the plane $t=0$, we
have that $\Hat{\gamma}'$, as a vector field on $\R^2$, is just a
translation of $\gamma'$ by the vector $(x_0,y_0)$.   Thus,
$\Hat{\gamma}$ is a seed curve for $\mS_0$.

\end{proof}

\medskip

We apply these results to a useful special case.

\medskip

\begin{cor}\label{C:seedparam}
Suppose that $\mS$ be a $H$-minimal entire graph over the
$(y,t)$-plane with empty characteristic locus, and that $\mS$ is not
itself a vertical plane. There exist a point $g_0 \in \mS$, an
interval $I\subset \R$, $\gamma \in C^3(I)$, $h_0 \in C^2(I)$, so
that a neighborhood $\mS_0$ of $g_0$ can be parameterized by
\eqref{para},\eqref{para2} for $s \in I$ and $r \in \R$.
\end{cor}

\begin{proof}[\textbf{Proof}]
Since we assume $\mS$ is a graph over the $(y,t)$-plane, there
exists $\psi\in C^2(\R^2)$ such that $\mS$ is described by
$x=\psi(y,t)$. Consider the defining function $\Psi(x,y,t) = x -
\psi(y,t)$. If $\Psi_t = \psi_t\equiv 0$, then we would have
$\psi(y,t) = f(y)$, and by the $H$-minimality of $\mS$, we would
conclude that $f(y) = \alpha y + \beta$, against the hypothesis that
$\mS$ is not a vertical plane. Therefore, there exists $g_0\in \mS$
such that $\psi_t(g_0) \neq 0$. The Implicit Function Theorem
implies the existence of a neighborhood of $g_0$ on $\mS$ which may
be written as a graph over the plane $t=0$ (with empty
characteristic locus). Applying Theorem \ref{T:thmA}, we obtain
intervals $I, \subset \R$, $\gamma \in C^3(I),h_0\in C^2(I)$, so
that a neighborhood of $g_0$ is parameterized by \eqref{para},
\eqref{para2} for $s\in I$, $r \in J$. To finish the proof, we need
to show that for every $s \in I$, we may extend the domain of $r$ to
the whole of $\R$. To see this, we note that for each $s_0 \in I$,
the curve
\[r\ \to\ \mathcal{L}_{(\gamma(s_0),h_0(s_0))}(r) =
(\gamma(s_0)+r\gamma'(s_0)^\perp,h(s_0,r))\ ,\] with $h(s,r)$ given
by \eqref{para2}, and $r \in J$, is a horizontal straight line
segment in $\HH$. Using \eqref{para2}, we see that the tangent
vector to this line is simply
\begin{equation}\label{tl}
\left(\gamma'(s_0)^\perp, - \frac{\gamma(s_0)\cdot
\gamma'(s_0)}{2}\right)\ =\ \gamma_2'(s_0) \; X_1(\mathcal F(s_0,r))
-\gamma_1'(s_0) \; X_2(\mathcal F(s_0,r))\ , \end{equation} which,
by Corollary \ref{C:seednux}, is precisely $(\nuX)^\perp$. Thus, by
the standard uniqueness theory for solutions to ordinary
differential equations, these line segments must coincide with
subsets of the horizontal line foliation of $\mS$ guaranteed by
Theorem \ref{foliation}.  As the entirety of these horizontal lines
are contained in $\mS$, we conclude that the parameterization given
by \eqref{para}, \eqref{para2} extends to $(s,r) \in I\times \R$.

\end{proof}

\medskip

In order to extract some crucial additional information from our
assumption that $\mS$ be a an entire graph with $\Sigma(\mS) =
\varnothing$, the following elementary lemma will be useful.

\medskip

\begin{lemma}\label{L:lines}
Let $g_1,g_2 \in \HH$, $v=(v_1,v_2,v_3),w=(w_1,w_2,w_3) \in \R^3$,
and consider the straight-lines $L_i = \{g_i+rv\mid r\in \R\}$,
$i=1,2$. If $\pi:\HH \ra \R^2$ denotes the projection to the
$(y,t)$-plane given by $\pi(z_1,z_2,t)=(0,z_2,t)$, then $\pi(L_1)
\cap \pi(L_2) = \varnothing$ if and only if  $\pi(g_1) \neq
\pi(g_2)$, and $v \times w$ is (Euclidean) perpendicular to
$(1,0,0)$.
\end{lemma}
\begin{proof}[\textbf{Proof}]
Suppose $\pi(L_1) \cap \pi(L_2) = \varnothing$, then it is obvious
that it must be $\pi(g_1) \neq \pi(g_2)$. Furthermore, we must also
have $\pi(v) \times \pi(w)=(v_2w_3-w_2v_3,0,0) =0$.  Thus, we
conclude that $(v_2,v_3)$ and $(w_2,w_3)$ are constant multiples of
one another. But then $v \times w$ takes the form
$(v_2w_3-w_2v_3,\star,\star)=(0,\star,\star)$, and we conclude that
$v \times w$ is (Euclidean) perpendicular to $(1,0,0)$, one
direction of the lemma. The opposite direction follows by simply
reversing the previous argument.

\end{proof}

\medskip

Next, we apply Lemma \ref{L:lines} in the case of the
parameterization given in Corollary \ref{C:seedparam}. The following
lemma plays a crucial role in the proof of Theorem
\ref{T:reductionI}.

\medskip

\begin{lemma}\label{L:constcurv}
Let $I\subset \R$, $\gamma \in C^3(I) ,h_0\in C^2(I)$, and consider
a portion $\mS_0$ of an $H$-minimal entire graph $\mS$ over the
$(y,t)$-plane having empty characteristic locus. Suppose that
$\mS_0$ be parameterized as in \eqref{para},\eqref{para2} for
$(s,r)\in I\times \R$. There exists a subinterval, $J \subset I$,
such that $\gamma(J)$ is either a straight line segment, or a
circular arc.
\end{lemma}

\begin{proof}[\textbf{Proof}]
By re-parameterizing $\gamma$, we may assume $0 \in I$.  As in the
proof of Corollary \ref{C:seedparam}, for every fixed $s\in I$, we
consider the horizontal straight lines contained in $\mS$, and
defined by
\begin{equation}\label{lineq}
\mathcal{L}_{(\gamma(s),h_0(s))}(r) =
\left(\gamma(s)+r\gamma'(s)^\perp,h_0(s)-\frac{r}{2}\gamma(s)\cdot\gamma'(s)\right)\
,\quad\quad\quad r\in \R\ . \end{equation}

For every $s \in I$ we consider the two lines
$L_1=\mathcal{L}_{(\gamma(0),h_0(0))}$,
$L_2=\mathcal{L}_{(\gamma(s),h_0(s))}$. The assumption that $\mS$ be
an entire graph over the $(y,t)$-plane implies, in particular, that
for every $s\in I$ we must have $\pi(L_1)\cap \pi(L_2)=
\varnothing$. We can thus use Lemma \ref{L:lines} to infer that the
directional vectors of $L_1$ and $L_2$,
$v=\left(\gamma'(0)^\perp,-\frac{1}{2}\gamma(0)
  \cdot \gamma'(0)\right)$ and
  $w=\left(\gamma'(s)^\perp,-\frac{1}{2}\gamma(s)
  \cdot \gamma'(s)\right)$, satisfy the condition that $v \times w$ be perpendicular to
  $(1,0,0)$.
 A simple computation now gives
\begin{equation}\label{eq0}
v \times w \cdot
(1,0,0)=\frac{1}{2}\bigg(\gamma_1'(0)(\gamma(s)\cdot
\gamma'(s))-\gamma_1'(s)(\gamma(0) \cdot \gamma'(0))\bigg)=0\ .
\end{equation}

We first discuss three special cases.  If $\gamma'_1(0)=0$ and
$\gamma(0)\cdot \gamma'(0) \neq 0$, then for \eqref{eq0} to be
satisfied we must have $\gamma_1'(s)=0$ for all $s$ in a
neighborhood $J$ of $s=0$.  Hence, $\gamma(J)$ is a line segment and
we have reached one of our conclusions. Similarly, if $\gamma(0)
\cdot \gamma'(0)=0$, but $\gamma_1'(0) \neq 0$, then \eqref{eq0} is
only satisfied if $\gamma(s) \cdot \gamma'(s) =0$ for all $s$ in a
neighborhood $J$ of $s=0$. In this case, we claim that $\gamma(J)$
is a circular arc. To see this, we note that $d/ds(|\gamma(s)|^2) =
2 \gamma(s) \cdot \gamma'(s)\equiv 0$ on $J$. Thus,
$|\gamma(s)|\equiv C >0$ on $J$, hence $\gamma(J)$ is a circular
arc, reaching the second of our conclusions. Last, if both
$\gamma(0) \cdot \gamma'(0)=0$ and $\gamma_1'(0) = 0$, then by the
fact that $|\gamma'(0)|=1$, we obtain that $|\gamma_2'(0)|=1$. and,
using \eqref{lineq}, we conclude that $\mS_0$ contains the straight
line
\[\mathcal{L}_{(\gamma(0),h_0(0))}(r)=(\gamma_1(0) \pm r,\gamma_2(0),h_0(0))\ ,\quad\quad\quad r\in \R\ .
\]

This clearly violates our assumption that $\mS_0$ is a graph over
the $(y,t)$-plane, as the whole line projects to the same point,
$(0,\gamma_2(0),h_0(0))$, in the $(y,t)$-plane, thus this case
cannot occur.

If we are not in any of these cases, then we must have that both
$\gamma_1'(0) \neq 0$ and $\gamma \cdot \gamma'(0) \neq 0$. By the
continuity of $\gamma$ and $\gamma'$, there exists a neighborhood of
$s=0$, $J$, such that $\gamma_1'(s) \neq 0$, and $\gamma(s) \cdot
\gamma'(s) \neq 0$ for any $s\in J$.  In this case we conclude that
\eqref{eq0} is equivalent to the existence, for every $s\in J$, of a
$C(s)\neq 0$ such that
\[
\begin{pmatrix}
\gamma_1'(0)
\\
\gamma_1'(s)
\end{pmatrix}\ =\ C(s)\ \begin{pmatrix} \gamma(0)\cdot \gamma'(0)
\\
\gamma(s) \cdot \gamma'(s)\end{pmatrix}\ .
\]

Equating the first entries we find that $C(s) \equiv C\neq 0$ for
$s\in J$, with
\begin{equation*}
C\ = \frac{\gamma_1'(0)}{\gamma(0) \cdot \gamma'(0)}\ .
\end{equation*}

We conclude that
\begin{equation}\label{eqa}
\gamma_1'(s)\ =\ C\ \gamma(s) \cdot \gamma'(s)\ ,\quad\quad\quad
s\in J\ .
\end{equation}

Recalling that $\frac{d}{ds} |\gamma(s)|^2 = 2 \gamma(s) \cdot
\gamma'(s)$, \eqref{eqa} implies
\begin{equation}\label{eqb}
|\gamma(s)|^2 = \frac{2}{C}\gamma_1(s)+C_0
\end{equation}
where $C_0$ is an integration constant.  We claim that this implies
that $\gamma(J)$ is a circular arc.  To see this, we left-translate
$\mS_0$ by $(-1/C,0,0)$. Lemma \ref{L:translateseed} shows that the
translated surface will have a seed curve given by
$\hat{\gamma}(s)=(\gamma_1(s)-1/C,\gamma_2(s))$.  Computing
$|\hat{\gamma}(s)|^2$ and using \eqref{eqb}, we have
\[|\hat{\gamma}(s)|^2=|\gamma(s)|^2-\frac{2}{C}\gamma_1(s)+\frac{1}{C^2} =
\frac{2}{C}\gamma_1(s)-\frac{2}{C}\gamma_1(s)+\frac{1}{C^2}+C_0 =
\frac{1}{C^2}+C_0\ .\]

Thus, $\hat{\gamma}(J)$ is a circular arc and hence, such is
$\gamma(J)$ as well.  This completes the proof.

\end{proof}

\medskip

With these preliminary computations in place, we are finally ready
to establish our main reduction result.

\medskip

\begin{proof}[\textbf{Proof of Theorem \ref{T:reductionI}}]
Let $\mS$ be an $H$-minimal entire graph over a plane $P$ with empty
characteristic locus and which is not itself a vertical plane. In
view of Lemma \ref{L:ytplane}, after possibly a left-translation and
a rotation about the $t$-axis, we may assume that $P$ is the plane
$x=0$, and that $\mS$ is given by $x=\psi(y,t)$ for some $\psi \in
C^2(\R^2)$. Corollary \ref{C:seedparam} guarantees that there exists
$g_0\in \mS$, an interval $I\subset \R$, a unit-speed $\gamma \in
C^3(I)$, $h_0 \in C^2(I)$, so that a neighborhood $\mS_0$ of $g_0$
can be parameterized by
\begin{equation}\label{seed0}
\mathcal F(s,r)\ =\ \left (\gamma(s) + r
\gamma'(s)^\perp,h_0(s)-\frac{r}{2}\gamma(s)
   \cdot \gamma'(s) \right )\ ,\quad\quad\quad (s,r)\in I\times \R\
   .
 \end{equation}

Since $\mS$ is a graph over the plane $x=0$, Lemma \ref{L:constcurv}
yields that $\gamma(J)$ is either a straight-line segment, or a
circular arc. We next show that assumption of empty characteristic
locus on $\mS$ rules out the possibility that $\gamma(J)$ be a
straight line segment. If, in fact, $\gamma$ were linear, then it
would have the form
\[\gamma(s)=(x_0+a_1s,y_0+a_2s)\ ,\]
with $a_1^2+a_2^2=1$. In this case, the signed curvature
$\kappa(s)\equiv 0$, and the angle function $W$, given by \eqref{W},
becomes
\[
W(r,s)\ =\ \left|h_0'(s) - r  + \frac{1}{2} (a_1y_0 - a_2
x_0)\right|\ ,\quad\quad\quad (s,r)\in J\times \R\ .
\]

Since it is clear that, for each fixed $s\in J$, there exists $r\in
\R$, $r=h_0'(s)+\frac{1}{2}(a_1y_0-a_2x_0)$, where $W(s,r) = 0$, we
would conclude that $\mS$ has a characteristic point at $g =
\mathcal F(s,r)$, against our hypothesis.

Therefore, $\gamma(J)$ must be a unit-speed circular arc, i.e.,
\begin{equation}\label{reduction}
\gamma(s) = (x_0+R\cos(s/R),y_0+R\sin(s/R))\ ,
\end{equation}
for some $(x_0,y_0)\in \R^2$, and $R>0$, and \eqref{seed0} becomes
\begin{equation}\label{seed}
\mathcal F(s,r)\ =\ \left
(x_0+(R+r)\cos(s/R),y_0+(R+r)\sin(s/R),h_0(s)+\frac{r}{2}(x_0\sin(s/R)-y_0\cos(s/R))\right
)\ .
\end{equation}

Consider the left-translated surface $\tilde \mS_0 =
(-x_0,-y_0,0)\circ \mS_0$ parameterized by
\begin{equation}\label{seed1}
\tilde{\mathcal F}(s,r) = (-x_0,-y_0,0)\circ \mathcal F(s,r)\ =\
\left((R+r)\cos \frac{s}{R},(R+r)\sin
\frac{s}{R},\overline{h}_0(s)\right)\ ,
\end{equation}
where
\[\tilde{h}_0(s) = h_0(s) + \frac{R}{2}\left(y_0\cos \frac{s}{R} - x_0\sin \frac{s}{R}\right)\ .\]

By Lemma \ref{L:translateplane}, we know that the $\tilde \mS =
(-x_0,-y_0,0)\circ \mS$ is a non-characteristic entire graph over a
plane $\tilde P$, having Euclidean normal $\tilde \bN^e =
(1,0,-y_0/2)$. Applying \eqref{ppar} to the parametrization
$\tilde{\mathcal F}(s,r)$ we obtain
\begin{equation}\label{ppar3}
\begin{cases}
\tilde p\  =\ -\ \sin \frac{s}{R}\ \left(\frac{(R + r)^2}{2R} -
\tilde h_0'(s)\right)\ ,
\\
\tilde q\  =\ \cos \frac{s}{R}\ \left(\frac{(R + r)^2}{2R} - \tilde
h_0'(s)\right)\ ,
\\
\tilde \omega\ =\ -\ \frac{R + r}{R}\ . \end{cases}
\end{equation}

Using \eqref{ppar3} we immediately recognize that the angle function
$\overline W$ for $\overline \mS_0$ is given by
\[
\overline W(s,r)\ =\ \left|\overline h_0'(s) -
\frac{(R+r)^2}{2R}\right|\ .
\]

Since $\overline W(s,r) \neq 0$ for any $(s,r)\in J \times \R$, we
conclude that we must either have
\[
\frac{(R+r)^2}{2R}\ >\ \tilde h_0'(s)\  ,\quad\text{or}\quad
\frac{(R+r)^2}{2R}\ <\ \tilde h_0'(s)
\]
for every $(s,r)\in J\times \R$. It is clear that for no fixed $s\in
J$ can the second inequality hold for every $r\in \R$, and therefore
we must have for every fixed $s\in J$ \[ \tilde h_0'(s)\ <\
\frac{(R+r)^2}{2R}\ ,\quad\quad\text{for every}\ r\in \R\ .
\]

This imposes that we must have the crucial property
\begin{equation}\label{injectivity}
\tilde h_0'(s)\ <\ 0\ ,\quad\quad\quad \text{for every}\quad s\in J\
.
\end{equation}

Having achieved this conclusion, consider the portion $\tilde \mS_0$
of $\tilde \mS$ given by \eqref{seed1}. The non-unit Euclidean
normal to $\tilde \mS_0$ is given by
\[
\tilde \bN^e\ =\ \tilde{\mathcal F}_s \times \ \tilde{\mathcal F}_r\
= \left(-\ \tilde h_0'(s)\ \sin \frac{s}{R}\ ,\ \tilde h_0'(s)\ \cos
\frac{s}{R}\ ,\ -\ \frac{R+r}{R}\right)\ ,
\]
where the cross product is taken with respect to the Euclidean inner
product in $\R^3$. The Euclidean Gauss map of $\tilde \mS_0$ is thus
given by
\[
\tilde \n^e\ =\  \frac{1}{\sqrt{\left(\frac{R+r}{R}\right)^2 +
\tilde h_0'(s)^2}}\ \left(-\ \tilde h_0'(s)\ \sin \frac{s}{R}\ ,\
\tilde h_0'(s)\ \cos \frac{s}{R}\ ,\ -\ \frac{R+r}{R}\right)\
\]

For any fixed $s \in J$, letting $r$ range over $\R$, we see that
the image of $\tilde \n^e$ in $S^2\subset \R^3$ describes a curve,
denoted $\Gamma$, which is an open arc of a great circle. We observe
that $\Gamma$ passes through the point $\frac{-
\hp}{|\hp|}(\sin(s),-\cos(s),0)= (\sin(s),-\cos(s),0)$, and that the
closure of $\Gamma$ contains the points $(0,0,\pm 1)$.  Since
$\tilde \mS_0$ is a graph over the plane $\tilde P$, we must have
that $\tilde \n^e(\tilde \mS_0) \subset S^2$ lies entirely to one
side of the plane determined by the vector
$\frac{1}{\sqrt{1+y_0^2/4}}(1,0,-y_0/2) \in S^2$. Now, if $y_0 \neq
0$, then the points $(0,0,\pm 1)\in \Gamma$ would lie on opposite
sides of this plane, thus reaching a contradiction. We conclude that
we must have $y_0=0$, and therefore $\tilde \mS_0$ is a graph over a
portion of the $(y,t)$-plane.

Since from \eqref{seed1} we see that $x/y = \cot(s/R)$, when $s\neq
0$, we fix an open sub-interval of $J$, $\tilde J = (a,b)$, with
either $0 <a <b<\pi$, or $-\pi < a<b<0$, and we consider the open
interval $I = \tilde h_0(\tilde J)$. We stress that, in view of
\eqref{injectivity} we know that $\tilde h_0^{-1} : I \to \tilde J$
exists, and, depending on the choice that we have made of $\tilde
J$, we have that either $0<\tilde h_0^{-1}(t) <\pi$, or $-\pi <
\tilde h_0^{-1}(t)<0$, for every $t\in I$. We may thus re-write
$\tilde \mS_0$ as
\begin{equation}\label{strip}
x\ =\ y\ \cot(\tilde{h}_0^{-1}(t))\ =\ y\ G(t)\  .
\end{equation}

We note that $G\in C^2(I)$ since $\tilde h_0^{-1}\in C^2(I)$, and
that moreover, thanks to \eqref{injectivity} we have
\[ G'(t)\ =\ -\
(1 + \cot^2(\tilde{h}_0^{-1}(t)))\frac{1}{\tilde h_0'(\tilde
h_0^{-1}(t))}\ >\ 0\ . \]

Furthermore, since $y=(R+r)\sin(s/R)$ for $s \in \tilde J$ and $r
\in\R$, we conclude that $y$ attains every real number. In addition,
we claim the map $(r,s)\in \R \times I \rightarrow
((R+r)\sin(s/R),\tilde{h}_0(s)) \in \R \times J$ is one to one. To
see this, we consider $(r_1,s_1),(r_2,s_2)$ so that
$((R+r_1)\sin(s_1/R),\tilde{h}_0(s_1))=((R+r_2)\sin(s_2/R),\tilde{h}_0(s_2))$.
The injectivity of $\tilde{h}_0$ implies that $s_1=s_2$, and so we
must have $r_1=r_2$ as well. We thus see that the parametrization
\eqref{strip} of $\tilde \mS$ is valid for $(y,t)\in \R \times I$.
 This completes the proof.

\end{proof}

\medskip

With Theorem \ref{T:reductionI} in hands, we can now establish our
main result of Bernstein type.

\medskip

\begin{proof}[\textbf{Proof of Theorem \ref{T:instability}}]
Let $\mS\subset \HH$ be a stable $H$-minimal entire graph, with
$\Sigma(\mS) =\varnothing$. Assume by contradiction that $\mS$ is
not a vertical plane. By Theorem \ref{T:reductionI}, after possibly
a left-translation and a rotation about the $t$-axis, the resulting
surface $\mS$ contains a strict graphical strip $\mS_0$. By Theorem
\ref{T:perminI} we know that $\mS_0$ is unstable, and therefore also
$\mS$ must be unstable, thus reaching a contradiction. We conclude
that $\mS$ must be a vertical plane.

\end{proof}

\vskip 0.6in


\section{\textbf{Obstruction to the higher-dimensional Bernstein
problem}}\label{S:bdg}

\vskip 0.2in

In this section we prove Theorem \ref{T:negative}. We begin with a
simple proposition, which is fact valid in any Carnot group.

\medskip

\begin{prop}\label{P:rigate}
Suppose that the hypersurface $\mathcal S\subset \Hn$ be a
\emph{vertical cylinder}, i.e., it can be represented in the form
\begin{equation}\label{vertical}
\mathcal S\ =\ \{g=(x,y,t)\in \Hn \mid \mathfrak h(x,y) = 0\}\ ,
\end{equation}
where $\mathfrak h \in C^2(\R^{2n})$, and there exist an open set
$\omega\subset \R^{2n}$ and $\alpha>0$ such that $|\nabla \mathfrak
h|\geq \alpha$ in $\omega$. Under these assumptions,
 the characteristic locus of $\mathcal S$ is empty, and the $H$-mean curvature of $\mathcal S$ is therefore globally defined and it is given by
\begin{equation}\label{Hvertical}
\mathcal H(x,y,t)\ =\ (2n-1)\ H(x,y)\ ,
\end{equation}
where $H(x,y)$ represents the Riemannian mean curvature of the
projection $\pi(\mathcal S)$ of $\mathcal S$ onto
$\R^{2n}\times\{0\}$. In particular, $\mathcal S$ is $H$-minimal
 if and only if $\pi(\mathcal S)$ is a classical minimal surface in
 $\R^{2n}$. Furthermore, the $H$-perimeter of $\mS$, $\sigma_H(\mS)$, is given by
\begin{equation}\label{eqper}
\sigma_H(\mS)\ =\ H^{2n}(\mS)\ ,
\end{equation}
where $H^{2n}$ is the standard $2n$-dimensional Hausdorff measure in
$\R^{2n+1}$, i.e., the surface measure.
\end{prop}

\begin{proof}[\textbf{Proof}]
Consider the defining function for $\mS$ given by $\phi(x,y,t) =
\mathfrak h(x,y)$. We first observe that for $i=1,...,n$,
\[
X_i\phi(x,y,t)\ =\ \frac{\partial \mathfrak h}{\partial x_i} \
,\quad\quad\quad X_{n+i}\phi(x,y,t)\ =\ \frac{\partial \mathfrak
h}{\partial y_i}\ ,
\]
hence, since $T\phi \equiv 0$, we conclude that
\begin{equation}\label{samegrad}
\nabh \phi\ =\ \nabla\phi\ =\ \nabla \mathfrak h\ ,
\end{equation}
which, thanks to the assumption $|\nabla\mathfrak h|\geq \alpha >0$,
proves in particular that $\Sigma(\mS) = \varnothing$, and that
$\nuX = \frac{\nabla \mathfrak h}{|\nabla \mathfrak h|} = \n$, where
$\n$ denotes the Riemannian unit normal of $\mS$. Thanks to
Proposition \ref{P:equalMC}, the $H$-mean curvature of $\mS$ is
given by
\[
\mathcal H\ =\ \sum_{i=1}^{2n}\ \di\ <\nuX,X_i>\ =\ \sum_{i=1}^{2n}
X_i \n_i\ =\ div\ \n\ =\ (2n-1)\ H\ ,
\]
where in the third and second to the last equality we have used
\eqref{samegrad}. This proves \eqref{Hvertical}. Finally,
\eqref{eqper} derives from \eqref{samegrad} and from \eqref{per}, or
equivalently \eqref{pm}.

\end{proof}

\medskip

\begin{proof}[\textbf{Proof of Theorem \ref{T:negative}}]
Consider the Heisenberg group $\Hn$, and denote by $N+1= 2n$ the
dimension of the horizontal layer $\R^{2n}\times \{0\}$. For
$(x,y)\in \R^{2n}$, we write $y=(y',y_n)$, with $y'\in \R^{n-1}$,
and denote by $\R^N = \R^n_x\times \R^{n-1}_{y'}$. By the
fundamental results in \cite{BDG}, given any $N\geq 8$ there exists
a non-affine $f\in C^\omega(\R^N)$ such that $\mS_0 = \{(x,y)\in
\R^{2n}\mid y_n = f(x,y')\}$ is an entire minimal graph. Clearly, if
we consider the defining function $\mathfrak h(x,y) = y_n - f(x,y')$
for $\mS_0$, then \[ |\nabla \mathfrak h(x,y')|\ =\ \sqrt{1 +
|\nabla_{x,y'} f(x,y')|^2}\ \geq\ 1\ ,\quad\quad\text{for
every}\quad (x,y')\in \R^N\ .
\]

Consider the vertical cylinder $\mS\subset \Hn$ such that $\pi(\mS)
= \mS_0$. Thanks to Proposition \ref{P:rigate}, $\mS$ is an
$H$-minimal entire graph, over the hyperplane $\{(x,y,t)\in \Hn\mid
y_n = 0\}$, with empty characteristic locus, and which is not a
vertical hyperplane. Using the fact that the unit vector field on
$\mS$
\[
(x,y,t)\ \to\ \nuX\ =\ \frac{1}{\sqrt{1 + |\nabla_{x,y'} f|^2}}\
\left\{\sum_{i=1}^n (-f_{x_i}) X_i\ +\ \sum_{i=1}^{n-1} (- f_{y_i})
X_{n+ i}\ + \ X_{2n}\right\}\ ,
\]
is independent of the $t$-variable, and moreover $div^H \nuX = div
\n = 0$, we can easily prove the stability of $\mS$ similarly to the
classical case, see \cite{CM}, pagg.1-4, and also \cite{BSV} for a
general discussion of sub-Riemannian calibrations in $\Hn$. Finally,
we observe that the condition $N\geq 8$, translates into $n\geq
9/2$, hence a counterexample to the Bernstein problem can be found
for any $n\geq 5$.

\end{proof}

\end{document}